\documentclass{amsart}
\usepackage{amssymb}
\usepackage{latexsym}
\usepackage{pstricks,pst-node}

\title{Alcoved Polytopes I}
\author{Thomas Lam and Alexander Postnikov}
\address{Department of Mathematics, M.I.T., Cambridge, MA, 02139}
\email{thomasl@math.mit.edu} \email{apost@math.mit.edu}
\thanks{A.P.\ was supported in part by National Science Foundation
grant DMS-0201494 and by Alfred P.\ Sloan Foundation research
fellowship}
\date{April~21, 2005; updated on February 2, 2006 and May 22, 2006}

\newtheorem{theorem}{Theorem}[section]
\newtheorem{proposition}[theorem]{Proposition}
\newtheorem{prop}[theorem]{Proposition}
\newtheorem{corollary}[theorem]{Corollary}
\newtheorem{cor}[theorem]{Corollary}
\newtheorem{definition}[theorem]{Definition}

\newtheorem{remark}[theorem]{Remark}
\newtheorem{lemma}[theorem]{Lemma}
\newtheorem{example}[theorem]{Example}

\def\tv{\bar v}
\def\tw{\bar w}
\def\R{\mathbb{R}}
\def\Z{\mathbb{Z}}
\def\C{\mathbb{C}}
\def\<{\left<}
\def\>{\right>}

\def\A{\mathcal{A}}

\def\H{\Delta}

\def\G{\Gamma}

\def\cdes{\mathrm{cdes}}

\def\des{\mathrm{des}}
\def\mod{\,\mathrm{mod}\,}

\def\I{{\mathcal{I}}}
\def\J{{\mathcal{J}}}
\def\ub{\underline{b}}
\def\uc{\underline{c}}

\newcommand{\floor}[1]{\lfloor #1 \rfloor }

\newcommand{\set}[1]{ \left\{ #1 \right\} }
\def\mm{{\mathcal M}}

\def\a{{\underline{a}}}

\def\N{{\mathbb N}}
\def\g{{\mathcal G}}
\def\A{{\mathcal A}}
\def\Z{{\mathbb Z}}
\def\b{{\underline{b}}}

\def\w{{\omega}}
\def\D{{\Delta}}
\def\L{{\mathcal L}}
\def\oo{{\mathcal O}}
\def\P{{\mathcal P}}

\def\cdes{\mathrm{cdes}}

\def\sort{\mathrm{sort}}

\psset{unit=0.7pt} \psset{arrowsize=4pt 1} \psset{linewidth=1pt}

\begin{document}

\begin{abstract}
The aim of this paper is to study alcoved polytopes, which are
polytopes arising from affine Coxeter arrangements. This class of
convex polytopes includes many classical polytopes, for example,
the hypersimplices. We compare two constructions of triangulations
of hypersimplices due to Stanley and Sturmfels and explain them in
terms of alcoved polytopes. We study triangulations of alcoved
polytopes, the adjacency graphs of these triangulations, and give
a combinatorial formula for volumes of these polytopes. In
particular, we study a class of matroid polytopes, which we call
the multi-hypersimplices.
\end{abstract}

\maketitle

\section{Introduction}

The affine Coxeter arrangement of an irreducible crystallographic
root system $\Phi\subset V\simeq \R^r$ is
obtained by taking all integer affine translations
$H_{\alpha,k} = \{x\in V\mid (\alpha,x) = k\}$, $\alpha\in\Phi$, $k\in\Z$,
of the hyperplanes perpendicular to the roots.  The regions of the affine
Coxeter arrangements are simplices called alcoves.  They are in a one-to-one
correspondence with elements of the associated affine Weyl group.
%If $\Phi$ is irreducible, then the alcoves are simplices
%congruent to each other.
We define an {\it alcoved polytope\/} $\P$ as a convex polytope that
is the union of several alcoves.   In other words, an alcoved polytope is the
intersection of some half-spaces bounded by the hyperplanes $H_{\alpha,k}$:
$$
\P= \{x\in V \mid b_\alpha\leq (\alpha,x)\leq c_\alpha,\ \alpha\in\Phi\},
$$
where $b_\alpha$ and $c_\alpha$ are some integer parameters.
These polytopes come naturally equipped with coherent
triangulations into alcoves.  Alcoved polytopes include many
interesting classes of polytopes: hypersimplices, order polytopes,
some special matroid polytopes, Fomin-Zelevinsky's generalized associahedra,
and many others.
This is the first of two papers about alcoved polytopes.  In this paper,
we concentrate on alcoved polytopes of the Lie type $A$ case
and on related combinatorial objects.  In~\cite{AP2}, we will
treat the general case of an arbitrary root system.

Hypersimplices are integer polytopes which appear in algebraic and geometric
contexts.  For example, they are moment polytopes for torus actions on
Grassmannians.  They are also weight polytopes of the fundamental
representations of the general linear group $GL_n$.
The $(k,n)$-th hypersimplex can be defined as the slice of the hypercube
$[0,1]^{n-1}$ located between the two hyperplanes $\sum x_i = k-1$
and $\sum x_i = k$.  It is well-know that the normalized volume of this
hypersimplex equals the Eulerian number $A_{k,n-1}$,
i.e., the number of permutations of size $n-1$ with $k-1$ descents.
Stanley~\cite{Sta1} explained this fact by
constructing a triangulation of the hypersimplex into $A_{k,n-1}$ unit
simplices.
Another construction of a triangulation of the hypersimplex was given by
Sturmfels~\cite{Stu}.  It naturally appears in the context of Gr\"obner bases.
These two constructions of triangulations are quite different.

In Section~\ref{sec:hypersimplex},
we compare these triangulations and show that they are actually
identical to each other and that they can be naturally described in terms of
alcoved polytopes.\footnote{C.~Haase~\cite{Haa} reported to us that he also
discovered this equivalence (unpublished).  G.~Ziegler~\cite{Zie} reported that
the alcove triangulation of the hypersimplex appeared as an example in his
1997/98 class on triangulations.} In Section~\ref{sec:triang_alcove},
we extend the
descriptions of this triangulation to general alcoved polytopes, and give a
formula for the volume of an alcoved polytope.  In
Sections~\ref{sec:matroid}--\ref{sec:multihyper}, we study in detail three
examples of alcoved polytopes: the matroid polytopes, the second hypersimplex
and the multi-hypersimplices.

In the second part~\cite{AP2} of this paper, we will extend the
hypersimplices to all Lie types and calculate their volumes.
We will prove a general theorem on volumes of alcoved polytopes.
We will give uniform generalizations of the descent and major index statistics, appropriate for our geometric
approach.

Some alcoved polytopes have also been studied from a more
algebraic perspective earlier; see~\cite{KKMS,BGT}.

%The current paper is the type $A$ portion of this general theory.

% We now give some notation.  Let $w \in S_n$ be a permutation.  We
% write $w_i$ for the image of $i$ under $w$.  A \emph{descent} of
% $w$ is an index $i \in \set{1,2,\ldots,n-1}$ such that $w_i >
% w_{i+1}$.
%
% The $(k,n)$-Eulerian number $A_{k,n}$ is the number of permutations
% $w \in S_n$ with $k-1$ descents.  The Eulerian polynomial $A_n(t)$
% is defined as
% \[
% A_n(t) = A_{1,n}t + A_{2,n}t^2 + \cdots + A_{n,n}t^{n}.
% \]
% for $n\geq 1$, and $A_0(q) = 1$.

%Throughout the paper, we will be considering integer polytopes
%$\P$ residing in $\R^n$ or in some hyperplane $H$ of $\R^n$.  In
%both cases the volume $\mathrm{Vol}(\P)$ will denote the
%\emph{normalized volume} of $\P$ with respect to the lattice
%$\Z^n$ or $\Z^n \cap H$ (see for example \cite{EC2}).  Thus for
%example the unit cube $[0,1]^n$ in $\R^n$ has volume $n!$.

%{\bf Acknowledgements.} Hi!
\section{Four triangulations of the hypersimplex}
\label{sec:hypersimplex}

Let us fix integers $0< k < n$.  Let $[n]:=\{1,\dots,n\}$ and
$\binom{[n]}{k}$
denote the collection of $k$-element subsets of $[n]$.
To each $k$-subset $I\in\binom{[n]}{k}$
we associate the $01$-vector $\epsilon_I =
(\epsilon_1,\dots,\epsilon_n)$ such that $\epsilon_i=1$, for $i\in I$;
and $\epsilon_i = 0$, for $i\not\in I$.

The {\it hypersimplex\/} $\H_{k,n}\subset\R^n$ is the convex polytope
defined as the convex hull of the points $\epsilon_I$,
for $I\in\binom{[n]}{k}$.
All these $\binom{n}{k}$ points are actually vertices
of the hypersimplex because they are obtained from each other
by permutations of the coordinates.
This $(n-1)$-dimensional polytope can also be defined as
$$
\H_{k,n} = \{(x_1,\dots,x_n)\mid 0\leq x_1,\dots,x_n\leq 1;\
x_1+\cdots + x_{n} =k\}.
$$
The hypersimplex is linearly equivalent to the
polytope $\tilde\Delta_{k,n}\subset\R^{n-1}$ given by
\[
\tilde\Delta_{k,n}=\{(x_1,\dots,x_{n-1})\mid 0\leq x_1,\dots,x_{n-1}\leq 1;\
k-1 \leq x_1 + \cdots + x_{n-1} \leq k\}.
\]
Indeed, the projection $p:(x_1,\dots,x_n)\mapsto (x_1,\dots,x_{n-1})$
sends $\Delta_{k,n}$ to $\tilde\Delta_{k,n}$.
The hypersimplex $\tilde\Delta_{k,n}$ can be thought of as
the region (slice) of the unit hypercube $[0,1]^{n-1}$
contained between the two hyperplanes $\sum x_i = k-1$ and $\sum x_i = k$.

Recall that a {\it descent\/} in a permutation $w\in S_n$ is an
index $i\in\{1,\dots,n-1\}$ such that $w(i)>w(i+1)$. Let $\des(w)$
denote the number of descents in $w$. The {\it Eulerian number\/}
$A_{k,n}$ is the number of permutations in $S_n$ with $\des(w)=k-1$
descents.
%The {\it Eulerian polynomial\/}
%$A_n(t)$
%is defined as
%$A_n(t) = A_{1,n}t + A_{2,n}t^2 + \cdots + A_{n,n}t^{n}$.
%for $n\geq 1$, and $A_0(q) = 1$.

Let us normalize the volume form in $\R^{n-1}$ so that the volume of
a unit simplex is 1 and, thus, the volume of a unit hypercube is
$(n-1)!$. It is a classical result, implicit in the work of 
Laplace~\cite[p.~257ff]{Lap},
that the normalized volume of the hypersimplex $\Delta_{k,n}$ equals
the Eulerian number $A_{k,n-1}$. One would like to present a
triangulation of $\Delta_{k,n}$ into $A_{k,n-1}$ unit simplices.
Such a triangulation into unit simplices is called a {\it unimodular
triangulation}.

In this section we define four triangulations of the hypersimplex
$\Delta_{k,n}$. One triangulation is due to Stanley~\cite{Sta1}, one
is due to Sturmfels~\cite{Stu}, one arises from the affine Coxeter
arrangement of type $A$ and the final one, which is new, we call the
circuit triangulation.  The main result of this section,
Theorem~\ref{thm:equaltriangulation}, says that these four
triangulations coincide.  In addition, we will describe the dual
graphs of these triangulations.

%As we explain in Section~\ref{sec:ranktwo}, such graphs have been
%studied previously for the case $k=2$; see~\cite{Dre,SY}.
%In the following sections we give four different constructions
%of such triangulations of the hypersimplex $\Delta_{k,n}$.
%One of these constructions is due to Stanley~\cite{Sta1}
%and another is due to Sturmfels~\cite{Stu}.

\subsection{Stanley's triangulation} \label{sec:Sta}

The hypercube $[0,1]^{n-1}\subset\R^{n-1}$ can be triangulated into
$(n-1)$-dimensional unit simplices $\nabla_w$ labelled by
permutations $w \in S_{n-1}$ given by
\[
\nabla_w = \set{(y_1,\ldots,y_{n-1}) \in [0,1]^{n-1} \mid 0 <
y_{w(1)} < y_{w(2)} < \cdots < y_{w(n-1)} < 1}.
\]

Stanley~\cite{Sta1} defined a transformation of the hypercube
$\psi:[0,1]^{n-1} \rightarrow [0,1]^{n-1}$ by
$\psi(x_1,\ldots,x_{n-1}) = (y_1,\ldots,y_{n-1})$, where
\[
y_i = (x_1 + x_2 + \cdots + x_i) - \lfloor x_1 + x_2 + \cdots +
x_i \rfloor.
\]
The notation $\lfloor x \rfloor$ denotes the integer part of $x$.
The map $\psi$ is piecewise-linear, bijective on the hypercube
(except for a subset of measure zero), and volume preserving.

Since the inverse map $\psi^{-1}$ is linear and injective when
restricted to the open simplices $\nabla_w$, it transforms the
triangulation of the hypercube given by $\nabla_w$'s into another
triangulation.

\begin{theorem}[Stanley~\cite{Sta1}]
The collection of simplices $\psi^{-1}(\nabla_w)$, $w\in S_{n-1}$,
gives a triangulation of the hypercube $[0,1]^{n-1}$ compatible
with the subdivision of the hypercube into hypersimplices.
The collection of the simplices $\psi^{-1}(\nabla_w)$,
where $w^{-1}$ varies over permutations in $S_{n-1}$
with $k-1$ descents, gives a triangulation of the $k$-th
hypersimplex $\tilde \Delta_{k,n}$.
Thus the normalized volume of $\tilde\Delta_{k,n}$
equals the Eulerian number $A_{k,n-1}$.
\label{thm:Sta}
\end{theorem}

\begin{proof}
Let $\psi(x_1,\dots,x_{n-1}) = (y_1,\dots,y_{n-1})\in\nabla_w$.
For $i=1,\dots,n-2$, we have,
$$
\floor{x_1+\cdots+x_{i+1}}=\left\{
\begin{array}{cl}
\floor{x_1+\cdots+x_i}&\textrm{if } y_i<y_{i+1};\\[.1in]
\floor{x_1+\cdots+x_i}+1&\textrm{if } y_i>y_{i+1}.\\[.1in]
\end{array}
\right.
$$
Thus $\floor{x_1+\cdots+x_{n-1}}=\des(w^{-1})$.
In other words, if $\des(w^{-1})=k-1$, then
$k-1\leq x_1+\cdots+x_{n-1}\leq k$, i.e., $(x_1,\dots,x_{n-1})\in
\tilde\Delta_{k,n}$.
\end{proof}

\subsection{Sturmfels' triangulation} \label{sec:Stu}

Let $S$ be a multiset of elements from $[n]$.  We define
$\sort(S)$ to be the unique non-decreasing sequence obtained by ordering
the elements of $S$.
%which has the same multiset of entries as $S$.
Let $I$ and $J$ be two $k$-element subsets of $[n]$,
and let $\sort(I \cup J) = (a_1,a_2,\dots, a_{2k})$.
Then we set $U(I,J) = \set{a_1,a_3,\dots,a_{2k-1}}$ and
$V(I,J) = \set{a_2,a_4,\dots,a_{2k}}$.
%If $I$ and $J$ are both $k$-subsets of $[n]$,
%then $U(I,J)$ and $V(I,J)$ are also
%$k$-subsets of $[n]$.
For example, for $I = \set{1,2,3,5}$, $J = \set{2,4,5,6}$, we have
$\sort(I\cup J) = (1,2,2,3,4,5,5,6)$,
$U(I,J) = \set{1,2,4,5}$, and $V(I,J) = \set{2,3,5,6}$.

We say that an ordered pair $(I,J)$ is {\it sorted} if $I = U(I,J)$ and
$J = V(I,J)$.  We call an ordered collection $\I=(I_1,\ldots,I_r)$
of $k$-subsets of $[n]$ \emph{sorted} if $(I_i,I_j)$ is sorted
for every $1 \leq i < j \leq r$.
Equivalently, if $I_l = \set{I_{l1}<\cdots<I_{lk}}$,
for $l=1,\dots,r$, then $\I$ is sorted if and only if
$I_{11} \leq I_{21} \leq \cdots \leq
I_{r1} \leq I_{12} \leq I_{22} \leq \cdots \leq I_{rk}$.
For such a collection $\I$, let $\nabla_\I$ denote the $(r-1)$-dimensional
simplex with the vertices $\epsilon_{I_1},\dots,\epsilon_{I_r}$.

\begin{theorem}[Sturmfels~\cite{Stu}]
The collection of simplices $\nabla_\I$, where $\I$ varies over all sorted
collections of $k$-element subsets in $[n]$,
is a simplicial complex that forms a triangulation
of the hypersimplex $\Delta_{k,n}$.
\end{theorem}

It follows that the maximal by inclusion sorted collections,
which correspond to the maximal simplices in the triangulation,
all have the same size $r=n$.

\begin{cor}
\label{cor:sorted_vol}
The normalized volume of the hypersimplex $\H_{k,n}$ is
equal to the number of maximal sorted collections of
$k$-subsets in $[n]$.
\end{cor}

This triangulation naturally appears in the context of Gr\"obner
bases.
%For the rest of this section assume reader's
%familiarity with Gr\"{o}bner bases and their relationship
%with coherent triangulations,
Let $k[x_I]$ be the polynomial ring in the $\binom{n}{k}$ variables
$x_I$ labelled by $k$-subsets $I\in\binom{[n]}{k}$. Define the map
$\phi: k[x_I] \rightarrow k[t_1,t_2,\ldots,t_n]$ by $x_I \mapsto
t_{i_1}t_{i_2}\cdots t_{i_k}$, for $I = \set{i_1,\ldots,i_k}$. The
kernel of this map is an ideal in $k[x_I]$ that we denote by
$\J_{k,n}$. Recall that a sufficiently generic height function on
the vertices $\epsilon_I$ of the hypersimplex $\Delta_{k,n}$ induces
a term order on monomials in $k[x_I]$ and defines a  Gr\"{o}bner
basis for the ideal $\J_{k,n}$.  On the other hand, such a height
function gives a coherent triangulation of $\Delta_{k,n}$. This
gives a correspondences between Gr\"obner bases and coherent
triangulations.  The initial ideal associated with a Gr\"obner basis
is square-free if and only if the corresponding triangulation is
unimodular. For more details on Gr\"obner bases, see
Appendix~\ref{sec:grobner_intro}.

%Under this correspondence, the radical of the
%initial ideal will equal to the Stanley-Reisner ideal of the
%triangulation.  It turns out that $\g_{k,n}$ is already
%radical, so the square-free standard monomials are in bijection
%with the simplices of the triangulation.

%The $\g_{k,n}$ for this ideal will be induced by some term order $\prec_{k,n}$
%which arises from some height function on the points $\epsilon_I$, giving a
%coherent triangulation of the convex hull $\conv\set{\epsilon_I}$ (.

\begin{theorem}[Sturmfels~\cite{Stu}]
\label{thm:Stu} The marked set of quadratic binomials
\[
\g_{k,n} = \set{\underline{x_I x_J} - x_{U(I,J)}x_{V(I,J)}\mid
I,J\in\binom{[n]}{k} },
\]
is a Gr\"{o}bner basis for $\J_{k,n}$ under some term order on
$k[x_I]$ such that the underlined term is the initial monomial.
The simplices of the corresponding triangulation are $\nabla_\I$,
where $\I$ varies over sorted collections of $k$-subsets of $[n]$.
Moreover, this triangulation is unimodular.
\end{theorem}

In Section~\ref{sec:sort_closed}, we state and prove a more
general statement.

%We will see in the next section that Stanley's triangulation and
%Sturmfels' triangulation are identical.

\subsection{Alcove triangulation}
\label{sec:alcove}

The {\it affine Coxeter arrangement\/} of type $A_{n-1}$ is the arrangement
of hyperplanes in $\R^{n-1}$ given by
$$
H_{ij}^l = \{(z_1,\dots,z_{n-1})\in\R^{n-1}\mid z_i-z_j = l\},\qquad\textrm{for }0\leq i<j\leq n-1,\ l\in\Z,
$$
where we assume that $z_0=0$.  It follows from the general theory
of affine Weyl groups, see~\cite{Hum}, that the hyperplanes
$H_{ij}^l$ subdivide $\R^{n-1}$ into unit simplices, called {\it
alcoves}.

We say that a polytope $\P$ in $\R^{n-1}$ is {\it alcoved\/} if
$\P$ is an intersection of some half-spaces bounded by the
hyperplanes $H_{ij}^l$. In other words, an alcoved polytope is a
polytope given by inequalities of the form $b_{ij} \leq z_j - z_i
\leq c_{ij}$, for some collection of integer parameters $b_{ij}$
and $c_{ij}$.  We will denote this alcoved polytope by $\P(b_{ij},
c_{ij})$.  If the parameters satisfy $b_{ij} = c_{ij} -1$, for all
$i,j$, then the corresponding polytope consists of a single alcove
(or is empty).  Each alcoved polytope comes naturally equipped
with the triangulation into alcoves.  Conversely, if $\P$ is a
convex polytope which is a union of alcoves, then $\P$ is an
alcoved polytope.
%An alcove is a non-empty alcoved
%polytope with parameters satisfying

Assume that $z_i = x_1+\cdots + x_i$, for $i=1,\dots,n-1$.
The hypersimplex $\tilde \Delta_{k,n}$ is given by
the following inequalities in the $z$-coordinates:
\begin{equation}
0\leq z_1-z_0,\dots,z_{n-1}-z_{n-2}\leq 1;\ k-1\leq z_{n-1}-z_0\leq k.
\label{eq:hypersimplex-z-coords}
\end{equation}
Thus the hypersimplex is an alcoved polytope.
Let us call its triangulation into alcoves the
{\it alcove triangulation}.

\subsection{Circuit triangulation} \label{sec:graphs}

Let $G_{k,n}$ be the directed graph on the vertices $\epsilon_I$,
$I\in\binom{[n]}{k}$, of the hypersimplex $\Delta_{k,n}$ defined as
follows. Let us regard the indices $i$ of a vector
$\epsilon=(\epsilon_1,\dots, \epsilon_n)$ as elements of $\Z/n\Z$.
Thus we assume that $\epsilon_{n+1}=\epsilon_1$. We connect a vertex
$\epsilon= (\epsilon_1,\dots,\epsilon_n)$ with a vertex $\epsilon'$
by an edge $\epsilon {\stackrel i \longrightarrow} \epsilon'$
labelled by $i\in[n]$ whenever $(\epsilon_i,\epsilon_{i+1})=(1,0)$
and the vector $\epsilon'$ is obtained from $\epsilon$ by switching
$\epsilon_i$ and $\epsilon_{i+1}$. In other words, each edge in the
graph $G_{k,n}$ is given by cyclically shifting a ``1'' in vector
$\epsilon$ one step to the right to the next adjacent place. It is
possible to perform such a shift if and only if the next place is
not occupied by another ``1''.

A circuit in the graph $G_{k,n}$ of minimal possible length is
given by a sequence of shifts of ``1''s so that the first ``1'' in
$\epsilon$ moves to the position of the second ``1'', the second ``1''
moves to the position of the third ``1'', and so on, finally, the last
``1'' cyclically moves to the position of the first ``1''. The length of
such a circuit is $n$.  We will call such circuits in $G_{k,n}$ is
{\it minimal}. Here is an example of a minimal circuit in
$G_{26}$:
$$
\begin{array}{ccccc}
(1,0,1,0,0,0) &\stackrel 3 \longrightarrow&
(1,0,0,1,0,0) &\stackrel 1 \longrightarrow&
(0,1,0,1,0,0) \\
\uparrow_6&&&& \downarrow_2\\
(0,0,1,0,0,1)
&\stackrel 5 \longleftarrow&
(0,0,1,0,1,0) &\stackrel 4 \longleftarrow&
(0,0,1,1,0,0)
\end{array}
$$
The sequence of labels of edges in a minimal circuit forms a
permutation $w=w_1\cdots w_n\in S_n$.
For example, the permutation corresponding to the above minimal
circuit is $w=3 1 2 4 5 6$.

If we do not specify the initial vertex in a minimal circuit, then
the permutation $w$ is defined modulo cyclic shifts $w_1\dots
w_n\sim w_{n}w_1\dots w_{n-1}$. By convention, we will pick the
representative $w$ of the class of permutations modulo cyclic
shifts such that $w_n=n$. This corresponds to picking the initial
point in a minimal circuit with lexicographically maximal
01-vector $\epsilon$. Indeed, if $\epsilon\stackrel i
\longrightarrow \epsilon'$ is an edge in $G_{k,n}$, then
$\epsilon>\epsilon'$ in the lexicographic order, for
$i=1,\dots,n-1$; and $\epsilon<\epsilon'$, for $i=n$.

\begin{lemma}
\label{lem:circuit-descents} A minimal circuit in the graph
$G_{k,n}$ is uniquely determined by the permutation $w$ modulo
cyclic shifts. A permutation $w \in S_n$ such that $w_n=n$
corresponds to a minimal circuit in the graph $G_{k,n}$ if and
only if the inverse permutation $w^{-1}$ has exactly $k - 1$
descents.
\end{lemma}
\begin{proof}
The inverse permutation $w^{-1}$ has a descent for each pair $(i <
j)$ such that $b = w_i = w_j + 1$.  Thus a ``1'' was moved from the
$b$-th position to the $(b+1)$-th position in $\epsilon$ before a
``1'' was moved from the $(b-1)$-th position to the $b$-th position.
(Here $\epsilon$ denotes the initial vertex of the corresponding
circuit.) This happens if and only if $\epsilon_b = 1$. Since
$\epsilon$ has $k$ ``1''s, this happens exactly $k$ times. But the
occurrence corresponding to $w^{-1}(n) = n > w^{-1}(1)$ is not
counted as a descent, so $w^{-1}$ has exactly $k-1$ descents.
Conversely, if $w^{-1}$ has $k-1$ descents, then we obtain a vector
$\epsilon$ with $k$ ``1''s, so that $w$ corresponds to a minimal
circuit containing $\epsilon$.
\end{proof}

For a permutation, $w=w_1\dots w_n\in S_n$, let $(w)$
denote the long cycle in $S_n$ given by $(w)=(w_1,\dots,w_n)$
in cycle notation.  Two permutations $u,w\in S_n$ are equivalent
modulo cyclic shifts if and only if $(u) = (w)$.
The reader should not confuse {\it circuits\/} in the graph $G_{k,n}$
with {\it cycles\/} in the symmetric group $S_n$.

Let $C_{k,n}$ denote the set of long cycles
$(w)=(w_1,\dots,w_{n-1},n)\in S_n$ such that $w^{-1}$ has exactly
$k - 1 $ descents. For $(w)\in C_{k,n}$, let $c_{(w)}$ be the
corresponding minimal circuit in the graph $G_{k,n}$, whose edges
are labelled by $w_1,\dots,w_n$. Lemma~\ref{lem:circuit-descents}
shows that that the map $(w)\mapsto c_{(w)}$ is one-to-one
correspondence between the set long cycles $C_{k,n}$ and the set of
minimal circuits in $G_{k,n}$.

Each minimal circuit $c_{(w)}$ in $G_{k,n}$ determines the simplex
$\Delta_{(w)}$ inside the hypersimplex $\H_{k,n}$ with the vertex
set $c_{(w)}$.

\begin{theorem} \label{thm:circuittriangulation}
 The collection of simplices
$\Delta_{(w)}$ corresponding to all minimal circuits in $G_{k,n}$
forms a triangulation of the hypersimplex $\H_{k,n}$.
\end{theorem}

Let us call this triangulation of the hypersimplex the {\it
circuit triangulation}.

\begin{theorem}  \label{thm:equaltriangulation} The following four triangulations of the hypersimplex
are identical:
Stanley's triangulation, Sturmfels' triangulation, the alcove triangulation,
and the circuit triangulation.
\end{theorem}

Let us prove Theorems~\ref{thm:circuittriangulation}
and~\ref{thm:equaltriangulation} together. Let $\Gamma_{k,n}$ be the
collection of (maximal) simplices of the triangulation of
Theorem~\ref{thm:equaltriangulation}.

\begin{proof}
The fact that Stanley's triangulation coincides with the alcove triangulation
follows directly from the definitions.  We leave this as an exercise
for the reader.

Let us show that the simplices $\Delta_{(w)}$ are exactly those
in Sturmfels' triangulation.  An ordered pair of subsets $I =
\set{i_1 <  \cdots < i_k}$ and $J = \set{j_1 < \cdots <
j_k}$ is sorted if and only if the interleaving condition $i_1
\leq j_1 \leq i_2 \leq j_2 \leq \cdots \leq j_k$ is satisfied.
When two vertices $\epsilon_I$ and $\epsilon_J$ belong to
the same minimal circuit, a ``1'' from $\epsilon_I$ is moved towards
the right in $\epsilon_J$ but never past the original position of
another ``1'' in $\epsilon_I$.
Thus the interleaving  $i_a \leq j_a \leq i_{a+1}$ condition is
satisfied, and similarly we obtain the other interleaving inequalities.
Conversely, the interleaving condition implies that
each sorted collection belongs to a minimal circuit in $G_{k,n}$.

Let us now show that the circuit triangulation coincides with
Stanley's triangulation.  Recall that the latter triangulation
occurs in the space $\R^{n-1}$. To be more precise, in order to
obtain Stanley's triangulation we need to apply the projection
$p:(x_1,\ldots,x_n) \mapsto (x_1,\ldots,x_{n-1})$ to the circuit
triangulation. Let us identify a permutation $w=w_1\cdots
w_{n-1}\in S_{n-1}$ with $k-1$ descents with the permutation
$w_1\cdots w_{n-1} n\in S_n$.

We claim that the projected simplex $p(\Delta_{(w)})$ is exactly
the simplex $\psi^{-1}(\nabla_w)$ in Stanley's triangulation.
Indeed, the map
$\psi^{-1}:(y_1,\dots,y_{n-1})\mapsto(x_1,\dots,x_{n-1})$
restricted to the simplex $\nabla_w = \set{0 < y_{w(1)} < \cdots <
y_{w(n-1)} <1}$, is given by $x_1 = y_1$ and
$$
x_{i+1} = \left\{
\begin{array}{cl}
y_{i+1}-y_i & \text{if } w^{-1}(i+1)>w^{-1}(i),\\
y_{i+1}-y_i+1 & \text{if }w^{-1}(i+1)<w^{-1}(i)
\end{array}
\right.
$$
for $i=1,\dots,n-2$. The vertices of the simplex $\nabla_w$ are
the points $v_0,\dots,v_{n-1}\in\R^{n-1}$ such that $v_r =
(y_1,\dots,y_{n-1})$ is given by $y_{w(1)} = \cdots = y_{w(r)} =
0$ and $y_{w(r+1)} = \cdots = y_{w(n-1)}=1$. The map $\psi^{-1}$
sends the vertex $v_0=(0,\dots,0)$ to the point
$(x_1,\dots,x_{n-1})$ such that $x_1=0$ and $x_{i+1} = 1$ if
$w^{-1}(i+1)<w^{-1}(i)$ and $x_{i+1} = 0$ if
$w^{-1}(i+1)>w^{-1}(i)$, for $i=1,\dots,n-2$.  The vertex $v_r$ is
obtained from $v_{r-1}$ by changing $y_{w(n-r)}$ from 0 to 1. Thus
$\psi^{-1}(v_{r})$ differs from $\psi^{-1}(v_{r-1})$ exactly in
the coordinates $x_{w(n-r)}$ and $x_{w(n-r)+1}$.  Here $x_n = k -
(x_1 + \ldots + x_{n-1})$.  In fact, going from
$\psi^{-1}(v_{r-1})$ to $\psi^{-1}(v_{r})$ we move a ``1'' from
$x_{w(n-r)+1}$ to $x_{w(n-r)}$.  Finally, moving from $v_{n-1}$ to
$v_0$ we are changing $x_1$ from $1$ to $0$.  Thus as we go from
$\psi^{-1}(v_r)$ to $\psi^{-1}(v_{r+1})$ we are traveling along
the edges of the graph $G_{k,n}$ in the reverse direction.  So the
vertices $\psi^{-1}(v_r)$ of the simplex $\psi^{-1}(\nabla_w)$ are
exactly the vertices of $p(\Delta_{(w)})$. This completes the
proof of the theorem.
\end{proof}

\begin{remark}
\label{rem:bijtheta} An explicit bijection $\theta$ between maximal
sorted collections of $k$-subsets of $[n]$ and permutations $w \in
S_n$ with $k-1$ descents satisfying $w_n=n$ can be constructed as
follows. Let $\I=(I_1,\ldots,I_n)$ be such a collection.  Every
number in $[n]$ must occur in $\bigcup_i I_i$. Set
$(a_1,\ldots,a_{kn}) = \sort(\bigcup_i I_i)$. Let $\alpha_k$ be such
that $a_{\alpha_k} = k$ and $a_{\alpha_k+1} = k+1$.  Then
$\theta(\I) = w_1w_2\cdots w_{n-1}n$, where $w_i \equiv \alpha_i
\pmod {n-1}$ with representatives taken from $[n-1]$. This bijection
is compatible with the correspondences in
Theorem~\ref{thm:equaltriangulation}. For example, $\alpha_i \mod
n-1$ tells us when a ``1'' is moved from $\epsilon_i$ to
$\epsilon_{i+1}$ where by convention a ``1'' is moved from
$\epsilon_{n}$ to $\epsilon_1$ in the last edge of a circuit.
\end{remark}

\subsection{Adjacency of maximal simplices in the hypersimplex}
\label{sec:adjacency}

Let us say that two simplices in a triangulation are {\it adjacent}
if they share a common facet.  Let us describe the adjacent simplices in
the triangulation $\G_{k,n}$, using first the construction of the
circuit triangulation.

\begin{theorem}
\label{thm:adjacency} Two simplices $\Delta_{(u)}$ and
$\Delta_{(w)}$ of $\G_{k,n}$ are adjacent if and only if there
exists $i=1,\dots,n$ such that $u_{i}-u_{i+1}\ne \pm 1\pmod n$ and
the cycle $(w)$ is obtained from $(u)$ by switching
$u_{i}$ with $u_{i+1}$, i.e., $(w) = (u_i,u_{i+1}) (u)
(u_i,u_{i+1})$. Here again we assume that $u_{n+1}=u_1$.
\end{theorem}

\begin{proof}
The two simplices $\Delta_{(u)}$ and $\Delta_{(w)}$ are adjacent
if and only if exactly one pair of their vertices differ.  This
means that the corresponding minimal circuits $c_{(u)}$ and
$c_{(w)}$ differ in exactly one place.  Let
$\epsilon'\stackrel{u_i} \longrightarrow \epsilon
\stackrel{u_{i+1}} \longrightarrow \epsilon''$ be three vertices
in order along the minimal cycle $c_{(u)}$.  Then we can obtain
another cycle $c_{(w)}$ from $c_{(u)}$ by changing only $\epsilon$
if and only if $u_{i}-u_{i+1}\ne \pm 1\pmod n$ so that
$\epsilon'\stackrel{ u_{i+1}} \longrightarrow \epsilon^*
\stackrel{u_{i}} \longrightarrow \epsilon''$ are valid edges. When
$u_{i}-u_{i+1} = \pm 1\pmod n$ we are either moving the same ``1''
twice or moving two adjacent ``1''s one after another.  In both
cases, the order of the shifts cannot be reversed, and so
$\epsilon$ cannot be replaced by another vertex.
\end{proof}
Alternatively, let $\I = (I_1,\ldots,I_n)$ be a sorted subset
corresponding to the maximal simplex $\nabla_\I$ of $\G_{k,n}$. Let
$t \in [n]$ and $I_t = \set{i_1,i_2,\ldots,i_k}$.  Then we can
replace $I_t$ in $\I$ by another $I'_t \in \binom{[n]}{k}$ to
obtain an adjacent maximal simplex $\nabla_{\I'}$ if and only if
the following holds.  We must have $I'_t = \set{i_1,\ldots, i'_a,
\ldots, i'_b,\ldots,i_k}$ for some $a \neq b \in [n]$ and $i'_a
\neq i'_b$, satisfying $i_a - i'_a = i'_b - i_b = \pm 1 \pmod n$
and also both $k$-subsets
$\set{i_1,\ldots,i'_a,\ldots,i_b,\ldots,i_k}$ and $\set{i_1 ,
\ldots, i_a, \ldots, i'_b, \ldots i_k}$ must lie in $\I$.  For
example, we may replace $\set{1,3,5,8}$ by $\set{1,2,6,8}$ if and
only if both $\set{1,2,5,8}$ and $\set{1,3,6,8}$ lie in $\I$.

\begin{figure}[ht]
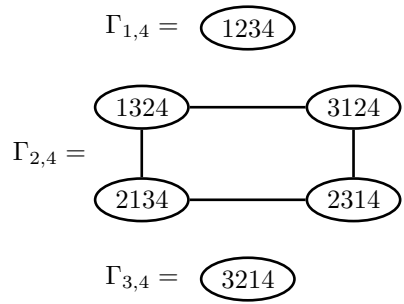

\pspicture(20,-15)(80,115)
\psset{unit=1pt}
\rput(10,80){$\G_{1,4}
= $} \rput(50,80){\ovalnode{A}{1234}} \rput(-25,32){$\G_{2,4} = $}
 \rput(10,50){\ovalnode{B}{1324}}
\rput(90,50){\ovalnode{C}{3124}} \rput(90,15){\ovalnode{D}{2314}}
\rput(10,15){\ovalnode{E}{2134}}

\rput(10,-15){$\G_{3,4}=$}\rput(50,-15){\ovalnode{F}{3214}}

\ncline{-}{B}{C} \ncline{-}{B}{E} \ncline{-}{C}{D}
\ncline{-}{D}{E}

\endpspicture
\caption{The graphs of the triangulations of $\H_{1,4}, \H_{2,4}$
and $\H_{3,4}$.} \label{fig:graphk4}
\end{figure}

\begin{figure}[ht]
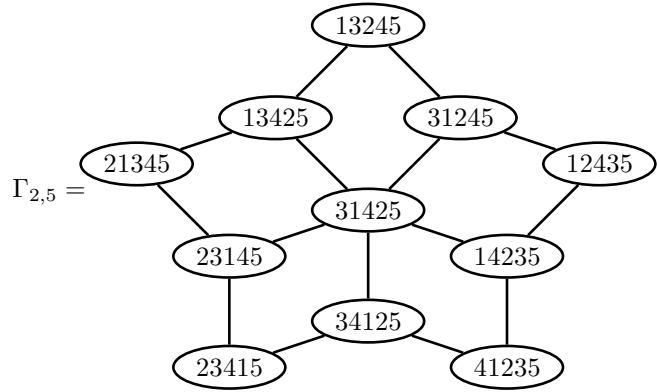

\pspicture(0 ,-10)(175,190)
\rput(100,20){\ovalnode{A}{34125}}
\rput(-70,90){$\G_{2,5} =$ } \rput(25,-5){\ovalnode{B}{23415}}
\rput(175,-5){\ovalnode{C}{41235}}
\rput(100,80){\ovalnode{D}{31425}}
\rput(25,55){\ovalnode{E}{23145}}
\rput(175,55){\ovalnode{F}{14235}}
\rput(-25,105){\ovalnode{G}{21345}}
\rput(225,105){\ovalnode{H}{12435}}
\rput(50,130){\ovalnode{I}{13425}}
\rput(150,130){\ovalnode{J}{31245}}
\rput(100,180){\ovalnode{K}{13245}} \ncline{-}{A}{B}
\ncline{-}{A}{C} \ncline{-}{A}{D} \ncline{-}{E}{D}
\ncline{-}{B}{E} \ncline{-}{F}{D} \ncline{-}{C}{F}
\ncline{-}{D}{I} \ncline{-}{D}{J} \ncline{-}{E}{G}
\ncline{-}{F}{H} \ncline{-}{G}{I} \ncline{-}{H}{J}
\ncline{-}{I}{K} \ncline{-}{J}{K}
\endpspicture
\caption{Graph of the triangulation of the hypersimplex
$\H_{2,5}$.} \label{fig:graph25}
\end{figure}

We can give $\G_{k,n}$ the structure of a graph by letting the
simplices be the vertices and letting an edge join two simplices
if the two simplices are adjacent.  Figures~\ref{fig:graphk4} and
\ref{fig:graph25} give examples of these graphs which we will also
denote as $\G_{k,n}$.

We will investigate degrees of vertices in these graphs for $k=2$ in
Section~\ref{sec:adjacency_second}.

\section{Triangulations and volumes of alcoved polytopes}
\label{sec:triang_alcove}
In this section, we generalize (Theorem~\ref{thm:alcove}) the
triangulations of Section~\ref{sec:hypersimplex} to all alcoved
polytopes (see Section~\ref{sec:alcove}).  As a consequence of these
descriptions of the alcove triangulation, we obtain a curious
formula (Theorem~\ref{thm:latticepoints}) expressing the volume of
an alcoved polytope as a sum of the number of lattice points in
certain other alcoved polytopes.  This formula has a root system
theoretic explanation which we give in~\cite{AP2}.  Finally, we
explain a construction of the dual graph of the alcove
triangulation, which we call the {\it alcove lattice}.

\subsection{Triangulations of alcoved polytopes}
\label{sec:alcoved}
%In the next few sections we study triangulations of general alcoved
%polytopes of type $A$. Most of the constructions here are
%generalized to other crystallographic root systems in~\cite{AP2}.

Let $\P = \P(b_{ij},c_{ij}) \subset \R^n$ be an alcoved polytope
which we realize in the $x$-coordinates. In other words, $\P$ is a
$(n-1)$-dimensional polytope lying in a hyperplane $x_1 + x_2 +
\cdots + x_n = k$ for some $k \in \Z$; and given by the inequalities
$b_{ij} \leq x_{i+1} + \cdots + x_j \leq c_{ij}$ for each pair
$(i,j)$ satisfying $0 \leq i < j \leq n-1$. By translating $\P$ by
$(m,\ldots,m)$ for some $m \in \Z$ to obtain an affinely equivalent
polytope, we can assume that all the coordinates of the points of
$\P$ are non-negative.  Let $Z_\P = \P \cap \Z^n \subset \N^n$
denote the set of integer points lying inside $\P$.

Let $G_\P$ be the directed graph defined as follows in analogy with
$G_{k,n}$ in Section~\ref{sec:graphs}.  The graph $G_\P$ has
vertices labelled by points $a \in Z_\P$.  Two vertices $a, b \in
Z_\P$ are connected by an edge $a \to b$ labelled $i$ if there
exists an index $i \in [1,n]$ such that $a + e_{i+1} - e_i = b$,
where $e_i, e_{i+1}$ are the coordinate vectors and $e_{n+1} :=
e_1$.  Let $C_\P$ denote the set of minimal circuits of $G_\P$,
which have length $n$.

For an integer vector $a = (a_1,a_2,\ldots,a_n) \in \N^n$ with
non-negative coordinates lying on $x_1 + x_2 + \cdots + x_n = k$, we
let $I_a$ denote the multiset of size $k$ of $\{1,2,\ldots,n\}$ with
$a_1$ 1's, $a_2$ 2's and so on.  If $I, J$ are multisets of size $k$
with elements from $\{1,2,\ldots,n\}$ then we can define $U(I,J)$
and $V(I,J)$ by sorting $I \cup J$ as in Section~\ref{sec:Stu}.
Similarly, we define the notions of sorted and sort-closed for
collections of multisets. In the following theorem, note that $Z_\P$
and $G_\P$ are defined without needing $\P$ to be alcoved.

The proof of the following theorem is exactly analogous to the
arguments of Theorem~\ref{thm:equaltriangulation}.

\begin{theorem}
\label{thm:alcove} Let $\P \subset \R^n$ be a $(n-1)$-dimensional
polytope lying in $x_1 + x_2 + \cdots + x_n = k$ so that all points
of $\P$ have non-negative coordinates.  Then the following are
equivalent:
\begin{enumerate}
\item
\label{it:1} $\P$ is an alcoved polytope.
\item
\label{it:3} The set $\mathcal{I} = \set{I_a \mid a \in Z_\P}$ is
sort-closed.  A triangulation of $\P$ consists of the maximal
simplices with vertices $\{a_1,a_2,\ldots,a_n\}$ for each sorted
collection $(I_{a_1},I_{a_2},\ldots,I_{a_n})$ where $I_{a_k} \in
\mathcal{I}$.
\item
\label{it:4} The set $C_\P$ of minimal circuits of $G_\P$ gives rise
to a triangulation of $\P$: if $C =(c^{(1)},c^{(2)},\ldots,c^{(n)})
\in C_\P$ then a maximal simplex is given by
$\mathrm{conv}(c^{(1)},c^{(2)},\ldots,c^{(n)})$.
\end{enumerate}
When these conditions hold, all the three above triangulations
agree.
\end{theorem}
%\begin{proof}
%The proof is exactly analogous to the arguments of
%Theorem~\ref{thm:equaltriangulation}.  The triangulation here is
%obtained by copying and translating the triangulations of
%Theorem~\ref{thm:equaltriangulation} by an integer vector so that
%the corresponding simplices cover the polytope $\P$.  For example,
%if $(c^{(1)},c^{(2)},\ldots,c^{(n)}) \in C_\P$ then there exists $c
%\in \Z^n$ so that $c^{(i)} - c$ is a 0-1 vector for each $i$. Thus
%$(c^{(1)}-c,c^{(2)}-c,\ldots,c^{(n)}-c)$ is a minimal circuit in
%some $G_{k,n}$ corresponding to some simplex $\Delta_{(w)}$; see
%Section~\ref{sec:graphs}.  The simplex $\Delta_{(w)} + c =
%\mathrm{conv}(c^{(1)},c^{(2)},\ldots,c^{(n)})$ will correspond to
%the circuit $(c^{(1)},c^{(2)},\ldots,c^{(n)})$.  To show that these
%simplices are disjoint and cover $\P$, we change to the
%$z$-coordinates of Section~\ref{sec:alcove} and see that we obtain
%the affine Coxeter arrangement.
%\end{proof}

An alcoved polytope also gives rise to a Gr\"{o}bner basis $\g_\P$
of the associated toric ideal $\J_\P$. The reader will be able to
write it down following Section~\ref{sec:sort_closed}.

Let us identify $w \in S_{n-1}$ with $w_1w_2\cdots w_{n-1} n \in
S_n$ as usual.  Recall that $\Delta_{(w)}$ denotes the simplex
$(\psi \circ p)^{-1}(\nabla_w)$ where we view $p$ as a map from $\{x
\in \R^n \mid x_1+\ldots+x_n = k\}$ to $\R^{n-1}$. For an alcoved
polytope $\P$, define the polytopes $\P_{(w)}$ by
\[
\P_{(w)} = \set{x \in \R^n \mid (\Delta_{(w)} + x) \subset \P}.
\]
Denote by $I(\P)$ the number $|\Z^n \cap \P| = \#Z_\P$ of lattice
points in $\P$.

\begin{theorem}
\label{thm:latticepoints} Each of the polytopes $\P_{(w)}$ is an
alcoved polytope.  The normalized volume of $\P$ is given by
\[
\mathrm{Vol}(\P) = \sum_w I(\P_{(w)})
\]
where the sum is over all permutations $w \in S_{n-1}$.
\end{theorem}
\begin{proof}
The alcoved triangulation is obtained by copying and translating the
triangulations of Theorem~\ref{thm:equaltriangulation} by an integer
vector so that the corresponding simplices cover the polytope $\P$.
For example, if $(c^{(1)},c^{(2)},\ldots,c^{(n)}) \in C_\P$ then
there exists $c \in \Z^n$ so that $c^{(i)} - c$ is a 0-1 vector for
each $i$. Thus $(c^{(1)}-c,c^{(2)}-c,\ldots,c^{(n)}-c)$ is a minimal
circuit in some $G_{k,n}$ corresponding to some simplex
$\Delta_{(w)}$; see Section~\ref{sec:graphs}.  The simplex
$\Delta_{(w)} + c = \mathrm{conv}(c^{(1)},c^{(2)},\ldots,c^{(n)})$
will correspond to the circuit $(c^{(1)},c^{(2)},\ldots,c^{(n)})$.
This proves the second statement of the theorem.

Let $\P$ be given by the inequalities $b_{ij} \leq x_{i+1} + \cdots
+ x_j \leq c_{ij}$ within the hyperplane $x_1+ x_2+\ldots+x_n = l$.
We check that $\P_{(w)}$ is an alcoved polytope.  In fact this
follows from the fact that $\Delta_{(w)}$ is itself an alcoved
polytope and is given by some inequalities $d_{ij} \leq x_{i+1} +
\cdots + x_j \leq f_{ij}$ and a hyperplane $x_1 + x_2 +\ldots+x_n =
k$, where we pick $d_{ij}$ and $f_{ij}$ so that all equalities are
achieved by some point in $\Delta_{(w)}$. Then $\P_{(w)}$ is the
intersection of the inequalities $b_{ij}-d_{ij} \leq x_{i+1} +
\cdots + x_j \leq c_{ij}-f_{ij}$ with the hyperplane $x_1+
x_2+\ldots+x_n = l-k$, which by definition is an alcoved polytope.
\end{proof}

For the hypersimplex, all the polytopes $\P_{(w)}$ are either empty
or a single point.

To conclude this section, we give one further interpretation of the
volumes of alcoved polytopes in terms of maps of the circle with
marked points.  Let $S^1$ be the unit circle and let $S^1_{(n)}$
denote a circle with $n$ distinct marked points
$p_0,p_1,\ldots,p_{n-1}$ arranged in clockwise order. Let $\P$ be an
alcoved polytope with parameters $b_{ij}$ and $c_{ij}$ as in
Section~\ref{sec:alcoved}. Let $M_\P$ denote the set of homotopy
classes of continuous maps $f: S^1_{(n)} \rightarrow S^{1}$
satisfying:
\begin{itemize}
\item
The map $f$ is always locally bijective and locally orientation
preserving.  Informally, this means that $f$ traces out $S^{1}$ in
the clockwise direction and never stops.
\item
The images of marked points are distinct.
\item
For each $0 \leq i < j \leq n-1$, The number $d$ of pre-images of
$f(p_i)$ under $f$ in the open interval $(p_i,p_j)$ satisfies
$b_{ij} \leq d < c_{ij}$.

\end{itemize}
Two maps $f$ and $g$ belong to the same homotopy class if and only
if they can be deformed into one another by a homotopy, in the usual
sense, while always satisfying the conditions above.  The following
proposition follows from the preceding discussion.

\begin{prop}
Let $\P$ be an alcoved polytope.  Then the simplices in the
triangulation $\G_\P$ are in bijection with the elements of $M_\P$.
\end{prop}
%\begin{proof}
%We think of $S^1$ as the interval $[0,2\pi]$ with the two endpoints
%identified.  For a point $x = (x_1,\ldots,x_n) \in \P$ in the
%interior of some alcove of $\P$, we obtain a map $f_x: S^1_{(n)}
%\rightarrow S^{1}$ by setting $f_x(p_0) = 0$ and between the points
%$p_i$ and $p_{i+1}$ we let $f_x$ be the unique injective map which
%travels a distance of $2\pi x_i$ at constant speed.

%If $x$ and $y$ lie in the same alcove of $\P$
%(Theorem~\ref{thm:alcove}) then $f_x$ and $f_y$ are easily seen to
%be homotopy equivalent, as the relative positions of the points
%$f_x(p_i)$ will not change. Conversely, picking an interior point in
%each alcove of $\P$ gives representatives of the homotopy classes of
%maps $f: S^1_{(n)} \rightarrow S^{1}$ satisfying the required
%conditions.
%\end{proof}

\subsection{Alcove lattice and alcoved polytopes} \label{sec:alcove
lattice}

Define the {\it alcove lattice} $\Lambda_n$ as the infinite graph
whose vertices correspond to alcoves (i.e., regions of type
$A_{n-1}$ affine Coxeter arrangement) and edges correspond to pairs
of adjacent alcoves. For example, $\Lambda_3$ is the infinite
hexagonal lattice. For an alcoved polytope $\P$, define its graph
$\Gamma_\P$ as the finite subgraph of $\Lambda_n$ formed by alcoves
in $\P$. For the graphs of Section~\ref{sec:graphs}, we have
$\Gamma_{k,n} = \Gamma_{\Delta_{k,n}}$.

According to~\cite[Sect.~14]{LP}, we have the following
combinatorial construction of the lattice $\Lambda_n$. Let
$[\lambda_0,\dots,\lambda_{n-1}]$ denote an element of
$\Z^n/(1,\dots,1)\Z$. In other words, we assume that
$[\lambda_0,\dots,\lambda_{n-1}]= [\lambda_0',\dots,\lambda_{n-1}']$
whenever the $\lambda_i'$ are obtained from the $\lambda_i$ by
adding the same integer. The vertices of $\Lambda_n$ can be
identified with the following subset of $\Z^n/(1,\dots,1)\Z$,
see~\cite{LP}:
$$
\Lambda_n = \{[\lambda_0,\dots,\lambda_{n-1}]\mid \textrm{the
integers } \lambda_0,\dots,\lambda_{n-1} \textrm{ have different
residues modulo } n\}.
$$
Two vertices $[\lambda_0,\dots,\lambda_{n-1}]$ and
$[\mu_0,\dots,\mu_{n-1}]$ of $\Lambda_n$ are connected by an edge
whenever there exists a pair $(i,j)$, $0\leq i\ne j\leq (n-1)$, such
that $\lambda_i+1\equiv \lambda_j\ \textrm{mod} \ n$ and
$(\mu_0,\dots,\mu_{n-1}) = (\lambda_0,\dots,\lambda_{n-1}) + e_i
-e_j$, where $e_i, e_j$ are the coordinate vectors in $\Z^n$. In
this construction, $\lambda_0,\dots,\lambda_{n-1}$ are the
$z$-coordinates of the central point of the associated alcove scaled
by the factor $n$, see~\cite{LP}. This construction immediately
implies the following description of the graph of the alcoved
polytope $\P(b_{ij},c_{ij})$, defined as in
Section~\ref{sec:alcove}.

\begin{proposition}  For an alcoved polytope $\P = \P(b_{ij},c_{ij})$,
its graph $\Gamma_\P$ is the induced subgraph of $\Lambda_n$ given
by the subset of vertices
$$
\left\{[\lambda_0,\dots,\lambda_{n-1}]\in\Lambda_n \mid n\cdot
b_{i,j}\leq \lambda_i - \lambda_j\leq n\cdot c_{i,j}  , \textrm{ for
} i,j\in[0,n-1] \right\}.
$$
\end{proposition}

The vertices of the graph $\Gamma_\P$ are in bijection with the
elements of $C_\P$ defined in the previous section.  Let $C =
(c^{(1)},c^{(2)},\ldots,c^{(n)}) \in C_\P$ be a minimal circuit in
the graph $\G_\P$.  The integer points
$\{c^{(1)},c^{(2)},\ldots,c^{(n)}\}$ are the vertices of an alcove
$A_C$, in the $x$-coordinates.  The vertex
$[\lambda_1,\lambda_2,\ldots,\lambda_n]$ of $\Lambda_n$ associated
to $A_C$ is given by $\lambda_i = \alpha_1 + \alpha_2 + \cdots +
\alpha_i$ where $\alpha := (\alpha_1,\ldots,\alpha_n)$ is given by
$\alpha = \sum_{i = 1}^n c^{(i)}$.  The point
$\frac{1}{n}[\lambda_1,\lambda_2,\ldots,\lambda_n] \in
\R^n/(1,\dots,1)\R$ is the central point of the alcove $A_C$ in the
$z$-coordinates.

\begin{example}
The $k$-th hypersimplex is given by
inequalities~(\ref{eq:hypersimplex-z-coords}).  Thus the vertex set
of the graph $\Gamma_{\Delta_{k,n}}$ is the subset of $\Lambda_n$
given by the inequalities $0\leq
\lambda_1-\lambda_0,\lambda_2-\lambda_1,\dots,\lambda_{n-1}-\lambda_{n-2}
\leq n$, and $(k-1)\cdot n\leq \lambda_{n-1}-\lambda_0\leq k\cdot
n$.
\end{example}

One can give an abstract characterization of subgraphs of
$\Lambda_n$ corresponding to alcoved polytopes. Let us say that an
induced subgraph $H$ of some graph $G$ is {\it convex} if, for any
pair of vertices $u,v$ in $H$, and any path $P$ in $G$ from $u$ to
$v$ of minimal possible length, all vertices of $P$ are in $H$.
In~\cite{AP2}, we will prove, in the more general context of an
arbitrary Weyl group, that an induced subgraph $\Gamma$ in
$\Lambda_n$ is the graph of some alcoved polytope if and only if
$\Gamma$ is convex.

\section{Matroid polytopes}
\label{sec:matroid}
Let $\mm$ be a collection of $k$-subsets of $[n]$.  The polytope
$\P_\mm$ is the convex hull in ${\mathbb R}^n$ of the points
$\set{\epsilon_I | \; I \in \mm}$ and is a subpolytope of the
hypersimplex $\H_{k,n}$.  In this section we classify the polytopes
$\P_\mm$ that are alcoved and in these cases give a combinatorial
interpretation for their volumes.

Our main interest lies in the case when $\mm$ is a matroid on the
set $[n]$.  In this case the polytopes $\P_\mm$ are known as {\it
matroid polytopes} -- we give explicit examples of matroid polytopes
which are alcoved.  Matroid polytopes have recently been studied
intensively, motivated by applications to tropical
geometry~\cite{FS,Spe}: the {\it amoeba} of a linear subspace $V
\subset \C^n$ is asymptotically described by its {\it Bergman fan},
and this fan is closely related to the normal fan of the matroid
polytope $\P_\mm$ of the matroid of $V$; see also~\cite{AK}.
Some of the results in this section have been obtained
earlier by Blum~\cite{Blu} in the context of Koszul algebras.
%where they were
%called base-sortable matroids.
%, from a more algebraic perspective.

\subsection{Sort-closed sets} \label{sec:sort_closed}
\begin{definition}
A collection $\mm$ of $k$-subsets of $[n]$ is \emph{sort-closed}
if for every two elements $I$ and $J$ in $\mm$, the subsets
$U(I,J)$ and $V(I,J)$ are both in $\mm$.
\end{definition}

A sorted subset of $\mm$ is a subset of the form
$\{I_1,\ldots,I_r\} \subset \mm$ such that $(I_1,\ldots,I_r)$ is a
sorted collection of $k$-subsets of $[n]$.

\begin{theorem}
\label{thm:matroid_tri} The triangulation $\G_{k,n}$ of the
hypersimplex induces a triangulation of the polytope $\P_\mm$ if
and only if $\mm$ is sort-closed. The normalized volume of
$\P_\mm$ is equal to the number of sorted subsets of $\mm$ of size
${\rm dim}(\P_\mm) + 1$.
\end{theorem}

The proof is analogous to that of Theorem~\ref{thm:Stu}.  We work in
the polynomial ring $k[x_I \mid I \in \mm]$. The ideal $\J_\mm$ is
the kernel of the ring homomorphism $\phi: k[x_I \mid I \in \mm]
\rightarrow k[t_1,t_2,\ldots,t_n]$ given by $x_I \mapsto
t_{i_1}t_{i_2}\cdots t_{i_k}$, for $I = \set{i_1,\ldots,i_k}$. The
following result is essentially equivalent to~\cite[Proposition 3.1]
{Blu}.

\begin{prop}
\label{prop:grobner_multi} Suppose that $\mm$ is sort-closed. Then
there is a term order, $\prec_\mm$, on $k[x_I |I \in \mm]$ so that
the reduced Gr\"{o}bner basis of $\J_\mm$ is given by the nonzero
marked binomials of the form
\begin{equation}
\label{bin} \set{\underline{x_I x_J} - x_{U(I,J)}x_{V(I,J)}}
\end{equation}

with the first monomial being the leading term.
\end{prop}
\begin{proof}[(Sketch of proof.)]
Since this proof is essentially the same as that of Sturmfels
in~\cite[Chapter 14]{Stu}, we will only sketch the argument; see
Appendix~\ref{sec:grobner_intro} for background.

We say a monomial $x_Ax_B\cdots x_V = x_{a_1\cdots a_k}x_{b_1\cdots
b_k} \cdots x_{v_1 \cdots v_k}$ is \emph{sorted} if the ordered
collection of sets $(A,B,\ldots,V)$ is sorted.  If a monomial is not
sorted, then there is a pair of adjacent variables $x_I x_J$ which
is unsorted.  Using the binomial $x_I x_J - x_{U(I,J)}x_{V(I,J)}$ we
can sort this pair.  We can sort a monomial modulo the ideal
generated by the marked binomials of~(\ref{bin}) in a finite number
of steps. Using~\cite[Theorem 3.12]{Stu}, we conclude that there is
a term order $\prec_\mm$ which selects the marked term for each
binomial of~(\ref{bin}).  Finally, one checks that the sorted
monomials are exactly the $\prec_\mm$-standard monomials.
\end{proof}

\begin{proof}[Proof of Theorem~\ref{thm:matroid_tri}]
The `if' direction follows from Proposition~\ref{prop:grobner_multi}
and Theorem~\ref{thm:regtri}.  We may
assume that $\prec_\mm$ arises from a weight vector since only
finitely many binomials are involved in (\ref{bin}).  For the
`only if' direction, suppose $\P = \P_\mm$ is a convex polytope
which is a union of simplices in $\G_{k,n}$. Since the
triangulation $\G_{k,n}$ is coherent, this triangulation $\G_\P$ is
also coherent.  We know already that all the faces of $\G_\P$ are
sorted collections of $k$-subsets of $[n]$ (if we identify a
simplex with its set of vertices).  By the correspondence of
Theorem~\ref{thm:regtri} and Proposition~\ref{prop:sqfree},
$\G_\P$ arises from some term-order $\prec_\mm$ which gives rise
to an initial ideal which is the Stanley-Reisner ideal of the
triangulation $\G_\P$.  Let $(I,J) \in \mm \times \mm$ not be
sorted, then $m_1 = x_Ix_J \in {\rm in}_{\prec_\mm}(\J_\mm)$.  As
in the proof of Proposition~\ref{prop:grobner_multi}, this means
there is another monomial $m_2$ so that $m_1 - m_2 \in \J_\mm$.
After a finite iteration of this argument, we see that $x_Ix_J =
x_Ax_B \ \mod \J_\mm $ for some $\prec_\mm$-standard monomial
$x_Ax_B$.  But this means that $x_Ax_B$ must be a sorted monomial
since it is an edge of an alcove.
  Thus $A = U(I,J)$ and $B= V(I,J)$ satisfy
$A, B \in \mm$ so that $\mm$ is sort-closed.
\end{proof}

%By Theorem~\ref{thm:equaltriangulation}, the set of vertices (viewed
%as a collection of $k$-subsets of $[n]$) of an alcoved subpolytope
%of the hypersimplex is sort-closed.  The converse follows from the
%description of an alcoved polytope as a convex union of alcoves; see
%Section~\ref{sec:alcove}.

\subsection{Sort-closed matroids} \label{sec:matroids} Let $k$ and
$n$ be positive integers satisfying $k \leq n$.  A (non-empty)
collection $\mm$ of $k$-subsets (called bases) of $[n]$ (the ground
set) is a \emph{matroid} if it satisfies the following axiom
(Exchange Axiom):
\begin{quote}
Let $I$ and $J$ be two bases of $\mm$.  Then for any $i \in I$
there exists $j \in J$ so that $(I - \set{i}) \cup \set{j}$ is a
base of $\mm$.
\end{quote}
The matroid $\mm$ is then said to be a rank $k$ matroid on $n$
elements.  If $I$ is a base of $\mm$ we write $I \in \mm$.  To a
matroid $\mm$ (of rank $k$ on $n$ elements) we associate the
\emph{matroid polytope} $\P_\mm$ as in
Section~\ref{sec:sort_closed}. We say that $\mm$ is {\it
sort-closed} if it is sort-closed as a collection of $k$-subsets.
Thus by Theorem~\ref{thm:matroid_tri}, the triangulation
$\G_{k,n}$ of the hypersimplex induces a triangulation of the
polytopes $\P_\mm$ for sort-closed matroids $\mm$.  Sort-closed
matroids were introduced by Blum~\cite{Blu}, who called them
base-sortable matroids.

We now describe two classes of matroids which are sort-closed. Let
$\Pi$ be a set partition of $[n]$ with parts
$\set{\pi_i}_{i=1}^{r}$ of sizes $|\pi_i|= a_i$, and $\ub =
(b_1,\ldots,b_r)$, $\uc = (c_1,\ldots,c_r)$ be two sequences of
non-negative integers. We will call the data $(\Pi,\ub,\uc)$ a
\emph{weighted set partition}. Define $\mm_{\Pi,\ub,\uc,k}$ to be
the collection of $k$-subsets $I$ of $[n]$ such that
\begin{equation} \label{eq:wsp} b_j \leq |I \cap \pi_j| \leq c_j
\end{equation} for all $j$.
\begin{lemma}
The collection of $k$-subsets $\mm_{\Pi,\ub,\uc,k}$ defined above
is a matroid.
\end{lemma}
\begin{proof}
Let $I$ and $J$ be two such subsets and $i \in I$, say $i \in
\pi_s$ for some $s$.  If $|I \cap \pi_k| = |J \cap \pi_k|$ for all
$k$ or if $|I \cap \pi_s| \leq |J \cap \pi_s|$, one can again find
some $j \in J \cap \pi_s - (I - \set{i}))$ to add to $I - \set{i}$
form a base. Otherwise there is some $t$ such that $|I \cap \pi_t|
< |J \cap \pi_t| \leq c_t$ in which case one can find some $j \in
(J \cap \pi_t)$ to add to $I - \set{i}$ without violating any of
the inequalities in (\ref{eq:wsp}). This verifies the exchange
axiom.
\end{proof}
This class of matroids is closed under duality.  The dual of
$\mm_{\Pi,\ub,\uc,k}$ is $\mm_{\Pi,\ub',\uc',n-k}$ where $b'_j =
|\pi_j|-c_j$ and $c'_j = |\pi_j| - b_j$.  We call the polytope
$\H_{\Pi,\ub,\uc,k}$ associated to the matroid $\mm_{\Pi,\ub,\uc,k}$
a {\it weighted multi-hypersimplex.} When $b_j=0$ and $c_j = 1$ for
all $j$ we will denote the matroid and polytope by $\mm_{\Pi,k}$ and
$\D_{\Pi,k}$ respectively, and we call the polytope $\D_{\Pi,k}$ a
\emph{multi-hypersimplex}.  Up to affine equivalence the polytope
$\D_{\Pi,k}$ depends only on the multiset $\set{a_i}_{i=1}^r$. The
polytope $\D_{\Pi,k}$ is the intersection of the hyperplane $x_1 +
\cdots + x_n = k$ with a product of simplices $\D_\Pi \simeq
\Delta_{a_1} \times \cdots \times \Delta_{a_r}$ just as the
hypersimplices are slices of cubes.  In the $z$-coordinates, this
polytope is determined by intersecting the hypersimplex
$\tilde{\H}_{k,n}$ with the inequalities
\begin{align*}
0 \leq  z_{a_1} - z_0 \leq 1;\;\; 0 \leq z_{a_1+a_2} - z_{a_1}
\leq 1; \ldots;  0 \leq z_{n} - z_{a_1+a_2+\cdots+a_{r-1}} \leq 1
\end{align*}
where we assume $z_n = k$.  A weighted multi-hypersimplex can be
viewed as a slice of a product of unions of hypersimplices.  In
particular, when $b_j = c_j -1$ for all $j$, the polytope
$\H_{\Pi,\ub, \uc,k}$ is the slice $x_1+\cdots+x_n=k$ of
\[
\H_{\Pi,\ub,\uc} \simeq \tilde{\H}_{c_1, a_1+1} \times
\tilde{\H}_{c_2,a_2+1} \times \cdots \times \tilde{\H}_{c_r,a_r+1}
\]
where the hypersimplex $\tilde{\H}_{c_i,a_i}$ (in the notation of
Section~\ref{sec:hypersimplex}) lives in the coordinates
$(x_{a_1+\cdots+a_{i-1}+1},\ldots,x_{a_1+\cdots+a_i})$.  Let
$\Pi(a_1,\ldots,a_r) $ denote the set partition
\[\set{\pi_1 = \set{1,\ldots,a_1},\ldots,
\pi_r=\set{a_1+\ldots+a_{r-1}+1,\ldots,a_1+\ldots+a_r = n}}.\]

\begin{prop}
The matroid $\mm_{\Pi,\ub,\uc,k}$ with $\Pi = \Pi(a_1,\ldots,a_r)$
and any $\ub, \uc \in \N^r$ is sort-closed.
\end{prop}
\begin{proof}
Let $I,J \in \mm_{\Pi,\ub,\uc,k}$.  Suppose to the contrary that
one of $U(I,J)$ or $V(I,J)$ were not a base.  Let
$(q_1,q_2,\ldots,q_{2k}) = \sort(I \cup J)$.  If $|U(I,J) \cap
\pi_s| > c_s$ or $|V(I,J) \cap \pi_s|> c_s$ then it must be the
case that for some $i$ the entries $q_i, q_{i+2}, \ldots,
q_{i+2c_s}$ belonged to the same part $\pi_s \in \Pi$. Then
$q_{i+1}, q_{i+3}, \ldots, q_{i+2c_s-1}$ belong to $\pi_s$ as well
since $q_{i+2k} \leq q_{i+2k+1} \leq q_{i+2k+2}$. This is
impossible as $I$ and $J$ were legitimate bases to begin with and
contain at most $c_s$ elements from $\pi_s$ each.  A similar
argument guarantees that $|U(I,J) \cap \pi_s| \geq b_s$ and
$|V(I,J) \cap \pi_s| \geq b_s$ for all $s$.
\end{proof}

A matroid $\mm$ is \emph{cyclically transversal} if it is a
transversal matroid specified by a set of (not necessarily disjoint)
cyclic intervals $\set{S_1,\ldots,S_k}$ of $[n]$.  Recall that the
bases of a transversal matroid are the $k$-element subsets $I =
\{i_1,\ldots,i_k\}$ of $[n]$ such that $i_s \in S_s$.

\begin{prop}[{\cite[Theorem 5.2 (without proof)]{Blu}}]
Let $\mm$ be a cyclically transversal matroid defined by the
subsets $\set{S_1,\ldots,S_k}$. Then $\mm$ is sort-closed.
\end{prop}
Since the proof is omitted in~\cite{Blu}, we give a simple direct
proof here.

\begin{proof}
Let $I$ be a $k$-element subset of $[n]$.  By the Hall marriage
theorem, $I$ is a base of $\mm$ if and only if \begin{equation}
\label{eq:hall} |I \cap \bigcup_{r\in R} S_r| \geq |R|
\end{equation} for every subset $R$ of $[k]$. Now let $I$ and $J$
be bases of $\mm$ and we now check (\ref{eq:hall}) for $U(I,J)$
and $V(I,J)$. Since each $S_i$ is a cyclic interval of $[n]$ it
suffices to consider the case where $\bigcup_{r\in R} S_r$ is
itself a cyclic interval $[a,b]$.  By hypothesis, the multiset $I
\cup J$ intersects $[a,b]$ in at least $2|R|$ elements.  Thus each
of $U(I,J)$ and $V(I,J)$ will intersect $[a,b]$ in at least $|R|$
elements.
\end{proof}
Let us describe $\P_\mm \subset \R^n$ for a cyclically transversal
matroid explicitly in terms of inequalities. It is given by the
hyperplane $x_1 + x_2 + \cdots + x_n = k$, the inequalities $0 \leq
x_i \leq 1$ together with the inequalities
\[
\sum_{s \in S_R} x_s \geq |R|
\]
with $S_R = \bigcup_{r\in R} S_r$ for every subset $R$ of $[k]$.

We end this section with the question: what other matroids are
sort-closed?

\section{The second hypersimplex}
There is a description of the triangulation of the second
hypersimplex, developed in~\cite{LST} and ~\cite[Chapter~9]{Stu}, in
terms of graphs known as thrackles.  In this section, we apply this
description to rank two matroids, and give a precise description of
the dual graph of the triangulation.
%and calculate the $f$-vectors of
%$\tilde{\G}_{2,n}$.

Triangulations of the second hypersimplex $\Delta_{2,n}$ arise in
the study of metrics on a finite set of points~\cite{Dre}, and
recently a thorough classification of the triangulations of
$\Delta(2,6)$ was performed in~\cite{SY}.  This classification of
triangulations is an important problem in phylogenetic
combinatorics.  Our study of the dual graphs of the triangulations
is partly motivated by this connection: the graph $\Gamma_{2,n}$ of
the triangulation is essentially what is known as the {\it tight
span} of the corresponding metric, and generalizes the phylogenetic
trees derived from the metric; see~\cite{Dre}.

\subsection{Thrackles and rank two matroids} \label{sec:ranktwo} When
the rank $k$ is equal to two (which we will assume throughout this
section), every matroid $\mm$ arises as $\mm_{\Pi,2}$ for some set
partition $\Pi$.  Throughout this section we will assume that $\Pi$
has the form $\Pi = \Pi(a_1,\ldots,a_r)$.
 Following~\cite[Chapter~9]{Stu} and \cite{LST}, we associate a graph
on $[n]$ to each maximal simplex of the matroid polytope
$\mm_{\Pi,2}$. The vertices are drawn on a circle so that they are
labelled clockwise in increasing order. Throughout this section, a
``graph on $[n]$'' will refer to such a configuration of the
vertices in the plane. Since the bases are two element subsets of
$[n]$, we may identify them with the edges.

\begin{lemma}
\label{lem:sortedge} Let $A = (a_1,a_2)$ and $B = (b_1,b_2)$ be
two bases. Then the pair $A, B$ is sorted if and only if the edges
$(a_1,a_2)$ and $(b_1,b_2)$ intersect (not necessarily in their
interior) when drawn on the circle.
\end{lemma}

\begin{proof}
Sorted implies that $a_1 \leq b_1 \leq a_2 \leq b_2$ which
immediately gives the lemma.
\end{proof}

Thus sorted subsets of $\binom{[n]}{2}$ correspond to graphs on
$[n]$ drawn on a circle, so that every pair of edges cross. These
graphs are known as {\it thrackles}.  Note that two edges sharing
a vertex are considered to cross.

\begin{prop}
\label{prop:volume_graph} Let $\Pi$ be a set partition and
$\D_{\Pi,2}$ have dimension $d$.  The maximal simplices in the
alcoved triangulation of $\D_{\Pi,2}$ are in one-to-one
correspondence with thrackles on $[n]$ with $d+1$ edges such that
all edges are bases.
\end{prop}
\begin{proof}
Follows immediately from Theorem~\ref{thm:matroid_tri} and
Lemma~\ref{lem:sortedge}.
\end{proof}
Without the condition on the number of edges in the thrackles of
Proposition~\ref{prop:volume_graph}, one would obtain graphs
corresponding to all simplices (not just the maximal ones) of the
triangulation.

\begin{figure}[ht]
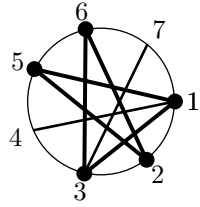

\pspicture(-50,-50)(50,50)
 \pscircle[linewidth=0.5pt](0,0){40}  \SpecialCoor
\psline[linewidth=1.5pt](40;0)(40;154)(40;308)(40;102)(40;257)(40;0)
\psline(40;51)(40;257) \psline(40;203)(40;0) \rput(50;0){$1$}
\rput(50;51){$7$}\rput(50;102){$6$}\rput(50;154){$5$}\rput(50;203){$4$}\rput(50;257){$3$}\rput(50;308){$2$}\rput(50;0){$1$}
\psdot[dotscale=1.5]
(40;0)\psdot[dotscale=1.5](40;154)\psdot[dotscale=1.5](40;308)\psdot[dotscale=1.5](40;102)
\psdot[dotscale=1.5](40;257)
\endpspicture
\caption{A thrackle $G(C)$.  The cycle $C$ has been drawn in
bold.} \label{fig:graphranktwo}
\end{figure}

When the dimension of $\D_{\Pi,2}$ is $n-1$, each thrackle $G$ is
determined by picking an odd-cycle $C$ such that all the edges cross
pairwise.  The remaining edges of $G$ join a vertex not on $C$ to
the unique `opposite' vertex lying on $C$ (so that the edge crosses
every edge of $C$);  see Figure~\ref{fig:graphranktwo}. We will call
the resulting thrackle $G(C)$. Let $C$ be a cycle, with pairwise
crossing edges, of length $2k+1$ with vertices $V(C) =
\set{v_1,v_2,\ldots, v_{2k+1}} \subset [n]$ labeled so that $v_1 <
v_2 < \cdots < v_{2k+1}$. Then the edges of $C$ are of the form
$(v_i,v_{k+i+1})$, where the indices are taken modulo $2k+1$. Thus
the condition that all the edges of $C$ are bases is equivalent to
$|V(C) \cap \pi_i| \leq k$ for all $i$. In fact this is enough to
guarantee that $G(C)$ corresponds to a valid maximal simplex of
$\D_{\Pi,2}$ -- that the remaining edges not on the cycle are bases
is implied.

Suppose $G$ arises from a sorted subset $(I_1,\ldots,I_r)$.  Let $w
= \theta(I_1,\ldots,I_r)$ where $\theta$ is the bijection of
Remark~\ref{rem:bijtheta}. The vertices $i$ not on the odd cycle of
$G$ are exactly the positions such that $w_{i} = w_{i-1} + 1$.

\begin{prop}
Let $a_1, \ldots, a_r$ be positive integers and $n = a_1 + \cdots +
a_r$. Then the $(n-1)$-dimensional volume of the second
multi-hypersimplex $\D_{\Pi(a_1,\ldots,a_r),2}$ is given by
\begin{align*}
\mathrm{Vol}(\P_{\Pi(a_1,\ldots,a_r),2}) &= \sum_{k=1}^\infty
\left(\sum_{c_1,..,c_r\leq k; \; c_1+...+c_r=2k+1}
\binom{a_1}{c_1} \cdots \binom{a_r}{c_r} \right) \\
&= 2^{n-1} - \sum_{i=1}^r \sum_{b,d\geq 0}\binom{a_i}{2b+d+1}
\binom{n-a_i}{d}
\end{align*}
\end{prop}

\begin{proof}
The first formula follows from enumerating odd subsets $S \subset
\{1,2,\ldots,n\}$ with size $2k+1$ satisfying $c_i := |S \cap
\pi_i| \leq k$ for all $i$.  The second formula comes from
counting the odd subsets $S' \subset \{1,2,\ldots,n\}$, where $|S'
\cap \pi_i| > |S'|/2$ for some $i \in \{1,\ldots,r\}$, and
subtracting them from all odd subsets of $\{1,\ldots,n\}$.
\end{proof}

One can also describe the simplices of the polytopes $\P_\mm$ for
higher rank matroids as hypergraphs $G$ satisfying the following
conditions:
\begin{enumerate}
\item
Every hyperedge $A \subset [n]$ of $G$ is a base of $\mm$.
\item
Let $A = \set{a_1 < \cdots < a_k}$, $B = \set{b_1 < \cdots < b_k}$
be a pair of hyperedges belonging to $G$.  Let $C_A$ be the cycle
on $[n]$ drawn with usual edges $(a_1,a_2), \ldots, (a_k,a_1)$ and
similarly for $C_B$.  Then each edge of $C_A$ must touch $C_B$ and
vice versa.
\end{enumerate}

\begin{remark}
Let $\tilde\Gamma_{2,n}$ be the
simplicial complex associated with the triangulation $\G_{2,n}$,
and let $f(\tilde\Gamma_{2,n})= \sum_{i=1}^d
f_i(\tilde\Gamma_{2,n})\, t^i$ be its $f$-polynomial,
where $f_i$ is the number of $i$-dimensional simplices in
this complex.
Using generating function techniques, one can deduce the following
expression for these polynomials:
\[\sum_{n\geq 2} f(\tilde\Gamma_{2,n},t)\, x^n =
\left[\frac{tq^2(1+q)(t^2 q^2 + t^2q -tq
+1)}{(1-tq)^2(1-2tq-tq^2)}\right]_{q \mapsto \frac{x}{1-x}}.\]
\end{remark}

\subsection{Adjacency of alcoves in the second hypersimplex}
\label{sec:adjacency_second} Let $\Delta \in \G_{k,n}$ be a
(maximal) simplex.  We say that $\Delta$ has \emph{degree} $d$ if it
is adjacent to $d$ other simplices.  We call $\Delta$ an
\emph{internal simplex} if none of its facets lies on the boundary
of $\H_{k,n}$. In this case $\Delta$ has maximal degree, namely $n$.

\begin{figure}
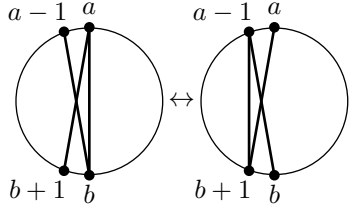

 \pspicture(-50,-50)(150,50)
 \pscircle[linewidth=0.5pt](0,0){40}
\pscircle[linewidth=0.5pt](100,0){40} \rput(0,50){$a$}
\rput(-28,48){$a-1$} \rput(-28,-48){$b+1$} \rput(0,-50){$b$}
 \rput(100,50){$a$}
\rput(72,48){$a-1$} \rput(72,-48){$b+1$} \rput(100,-50){$b$}
\rput(50,0){$\leftrightarrow$} \SpecialCoor \psline(40;90)(40;270)
\psline(40;90)(40;250) \psline(40;270)(40;110) \psdot(40;90)
\psdot(40;250)\psdot(40;270)\psdot(40;110)

\psset{origin={-100;0}} \psline(40;250)(40;110)
\psline(40;90)(40;250) \psline(40;270)(40;110) \psdot(40;90)
\psdot(40;250)\psdot(40;270)\psdot(40;110)
\endpspicture
\caption{The two simplices in $\G_{2,n}$ corresponding to two
thrackles $G$ and $G'$ are adjacent if $G$ and $G'$ are related by
the above move.}\label{fig:graph_move}
\end{figure}

\begin{proposition}
\label{prop:graph_move} The two simplices of $\G_{2,n}$
corresponding to two thrackles $G$ and $G'$ (via the
correspondence of Proposition~\ref{prop:volume_graph}) are
adjacent if and only if there are four distinct vertices labelled
$a, a+1, b, b+1 \ \mod n$ such that $G$ contains the edges
$(a,b),(a-1,b),(b+1,a)$ and $G'$ is obtained from $G$ by changing
the edge $(a,b)$ to $(a-1,b+1)$;  see Figure~\ref{fig:graph_move}.
\end{proposition}
\begin{proof}
The proposition follows immediately from Proposition~\ref{prop:volume_graph}
and Theorem~\ref{thm:adjacency} applied to
the case $k=2$ (more precisely, the comments after the proof of
the theorem).
\end{proof}

In fact the proposition is true also for smaller dimensional faces
of the simplices in $\G_{2,n}$.  Let $C$ be an odd cycle such that
all edges cross and let $|C|$ denote its length.

\begin{theorem}
\label{thm:degree_of_simplex} Let $\Delta$ be the simplex of
$\G_{2,n}$ corresponding to the thrackle $G(C)$.  If $G(C)$ is a
triangle then $\Delta = \H_{2,3}$ is the unique simplex in
$\G_{2,3}$ (and has degree 0).  If $|C| = 3$ and $G(C)$ has two
vertices of degree two then $\Delta$ has degree two.  Otherwise,
$\Delta$ has degree $|C|$.
\end{theorem}
\begin{proof}
The Theorem follows from Proposition~\ref{prop:graph_move}. Indeed
we can perform the move shown in Figure~\ref{fig:graph_move} to an
edge $(a,b)$ if it joins two vertices $a,b$ each of degree at
least two, with the exception of the case where $a$ and $b$ both
have degree exactly two and are joined to the same vertex $c$.
When the edge $(a,b)$ can be replaced, the change is necessarily
unique.  The case where $a$ and $b$ both have degree two and are
joined to the same vertex $c$ occurs only when $C$ is a
three-cycle.  In all other cases, every edge of $C$ can be
replaced, and so $\Delta$ has degree $|C|$.
\end{proof}

The following corollary is immediate from Theorem~\ref{thm:degree_of_simplex}.

\begin{corollary}
For $d > 1$, the second hypersimplex $\G_{2,n}$ has $\binom{n}{2d+1}$
simplices with degree $2d+1$.  No simplex has even degree
greater than or equal to 4.  In particular, if $n$ is odd then
$\G_{2,n}$ contains a unique internal simplex. If $n$ is even,
then $\G_{2,n}$ has no internal simplices.
\end{corollary}

\section{Multi-hypersimplices and multi-Eulerian polynomials}
\label{sec:multihyper}
%The hypersimplex is a slice of a cube, or alternatively, a slice of
%the product $[0,1] \times [0,1] \times \cdots \times [0,1]$ of
%intervals.

In this section we investigate the volumes of the {\it
multi-hypersimplices}, defined in Section~\ref{sec:matroids}. They
are slices of the product $\Delta_{a_1} \times \Delta_{a_2}\times
\cdots \times \Delta_{a_r}$ of simplices or more generally of
hypersimplices.  We define the multi-Eulerian numbers to be the
volumes of these polytopes.  They are generalizations of the usual
Eulerian numbers. These numbers, like the usual Eulerian numbers,
satisfy a number of interesting enumerative identities.  In the
first non-trivial case (Proposition~\ref{prop:multiEulerian}), we
determine the multi-Eulerian numbers explicitly in terms of Eulerian
numbers.

\subsection{Descent-restricted permutations and alcoved polytopes}
\label{sec:descent}

Let us consider an alcoved polytope $\P$ (as in
Section~\ref{sec:alcove}) which lies within a hypersimplex
$\H_{k,n}$.  In the $x$-coordinates, $\P \subset \R^n$ is defined by
the hyperplane $x_1 + x_2 + \cdots + x_n = k$, the inequalities $0
\leq x_i \leq 1$ together with inequalities of the form
\[
b_{ij} \leq x_{i+1} + \cdots + x_j \leq c_{ij}
\]
for integer parameters $b_{ij}$ and $c_{ij}$ for each pair $(i,j)$
satisfying $0 \leq i < j \leq n-1$.

Let $W_\P \subset S_{n-1}$ be the set of permutations $w =
w_1w_2\cdots w_{n-1} \in S_{n-1}$ satisfying the following
conditions:
\begin{enumerate}
\item
$w$ has $k-1$ descents.
\item
The sequence $w_{i} \cdots w_{j}$ has at least $b_{ij}$ descents.
Furthermore, if $w_{i} \cdots w_{j}$ has exactly $b_{ij}$ descents,
then $w_{i} < w_{j}$.
\item
The sequence $w_{i} \cdots w_{j}$ has at most $c_{ij}$ descents.
Furthermore, if $w_{i} \cdots w_{j}$ has exactly $c_{ij}$ descents,
then we must have that $w_{i} > w_{j}$.
\end{enumerate}
In the above conditions we assume that $w_0 = 0$.

Let $p: \R^n \rightarrow \R^{n-1}$ denote the projection as in
Section~\ref{sec:hypersimplex}. We may apply Stanley's
piecewise-linear map $\psi$ to $p(\P)$.  Using
Theorem~\ref{thm:equaltriangulation}, we see that a unimodular
triangulation of $p(\P)$ is given by the set of simplices
$\psi^{-1}(\nabla_{w^{-1}})$ as $w = w_1w_2\cdots w_{n-1}$ varies
over permutations in $W_\P$.  As a corollary we obtain the following
result.

\begin{proposition}
\label{prop:descent} The volume of $\P$ is equal to $|W_\P|$.
\end{proposition}

%In particular, $\G_{k,n}$ induces a triangulation of $\P$.

%\begin{proof}
%By Theorem~\ref{thm:equaltriangulation}, the polytope $p(\P)$ is the
%union of a set of simplices of the form $\psi^{-1}(\nabla_{w^{-1}})$
%for $w \in S_{n-1}$.  We check that the three conditions defining
%$W_\P$ hold.  The first condition follows as in the proof of
%Theorem~\ref{thm:Sta}.  As before, write $z_i = x_1 + \cdots + x_i$
%and $y_i = z_i - \floor{z_i}$.  We assume $z_0 = 0$ in addition. The
%inequality $0 \leq x_i \leq 1$ gives $0 \leq z_{i+1}-z_i \leq 1$.
%For a generic point, the number of descents in the sequence
%$(y_i,y_{i+1},\ldots,y_j)$, for $0 \leq i < j \leq n-1$, is equal to
%the number of indices $l \in [i,j-1]$ such that $z_l < a < z_{l+1}$
%for a whole number $a$.  If $y \in \nabla_{w^{-1}}$ then this is
%also the number of descents in the sequence $w_iw_{i+1}\cdots w_j$.
%Since we have $b_{ij} \leq z_j - z_i \leq c_{ij}$ the latter two
%conditions defining $W_\P$ follow. Conversely, that $w \in W_\P$
%implies $\psi^{-1}(\nabla_{w^{-1}}) \subset p(\P)$ follows in the
%same manner.
%\end{proof}
\subsection{Multi-Eulerian polynomials} \label{sec:multi} Let $\Pi$
be a fixed set partition with parts of sizes $\set{a_i}_{i=1}^{r}$
of total size $n = \sum_i a_i$ and let $\ub = (b_1,\ldots,b_r), \uc
= (c_1,\ldots,c_r) \in \N^r$. Define the \emph{weighted
multi-Eulerian number} $A_{\Pi,\ub,\uc,k} =
\mathrm{Vol}(\D_{\Pi,\ub,\uc,k})$ as the normalized
$(n-1)$-dimensional volume of the corresponding weighted
multi-hypersimplex.  We will consider polytopes of smaller dimension
to have volume 0 for what follows.  Now define the \emph{weighted
multi-Eulerian polynomial} $A_{\Pi,\ub,\uc}(t)$ as
\[
A_{\Pi,\ub,\uc}(t) = \sum_{k=1}^{r} A_{\Pi,\ub,\uc,k}t^{k}.
\]
Note that when $\Pi$ is the set partition with all $a_i = 1$ and
$b_i = 0$ and $c_i = 1$ for all $i$, then $A_{\Pi,\ub,\uc}(t)$
reduces to the usual Eulerian polynomial $A_{n-1}(t) = \sum_k
A_{k,n-1} t^k$. If $\Pi$ is a set partition of $[n]$, we denote by
$\Pi^*$ the set partition of $[n+1]$ with an additional part of
size one containing $n+1$. If $\underline{d} = (d_1,\ldots,d_r)
\in \Z^r$, then we denote by $\underline{d}^* \in \Z^{r+1}$ to be
the integer vector with an additional coordinate $d_{r+1} = 0$ and
$\underline{d}' \in \Z^{r+1}$ similarly with $d_{r+1} = 1$. Then
$\D_{\Pi^*,\ub^*,\uc',k}$ is affinely equivalent to (and has the
same normalized volume as) the intersection of $k-1 \leq
x_1+\cdots + x_n \leq k$ with the product of unions of
hypersimplices
\[
\left(\bigcup_{i = b_1+1}^{c_1}\H_{i,a_1+1}\right) \times \cdots \times \left(\bigcup_{i = b_r+1}^{c_r}\H_{i,a_r+1}\right).
\]
Thus
\[
A_{\Pi^*,\ub^*,\uc'}(1) =
\binom{n}{a_1,a_2,\ldots,a_r}\left(\sum_{i=b_1+1}^{c_1}A_{i,a_1}\right)\cdots\left(\sum_{i=b_r+1}^{c_r}A_{i,a_r}\right).
\]

In particular if $b_j = 0$ and $c_j = 1$ then this value is simply
a multinomial coefficient.  For this special case, we will omit
$\ub$ and $\uc$ in the notation and omit the prefix
\emph{weighted} from the names. We write out the combinatorial
interpretation for $A_{\Pi,k}$ explicitly.

\begin{prop}
\label{prop:multiperm} Let $k$ be a positive integer and $\Pi$ be
a set partition of $[n]$ with parts of sizes $\set{a_i}_{i=1}^{r}$
as before.  The following quantities are equal:
\begin{enumerate}
\item the multi-Eulerian number $A_{\Pi,k}$,
\item
\label{it:des} the number of permutations of $w \in S_{n-1}$ with
$k-1$ descents such that the substring $w_{a_1}w_{a_1+a_2} \cdots
w_{a_1+a_2+\cdots+a_{r-1}}$ has $k-1$ descents,
\item the number of sorted subsets $(I_1,\ldots,I_n)$ of $\mm_\Pi$ of size $n$.
\end{enumerate}
\end{prop}
\begin{proof}
As described in Section~\ref{sec:matroids} we may consider
multi-hypersimplices subpolytopes of the hypersimplex so the
proposition follows from Theorems~\ref{thm:equaltriangulation},
\ref{thm:matroid_tri} and Proposition~\ref{prop:descent}.
\end{proof}

Note that with $\Pi = \Pi(a_1,a_2,\ldots,a_r)$, then $A_{\Pi,k}$
is a function symmetric in the inputs $\{a_i\}$.  Let $\Pi=
\Pi(a,1^{n-a})$. By
Proposition~\ref{prop:multiperm}(\ref{it:des}), $A_{\Pi}(t)$ is
the generating function by descents for permutations satisfying
$w_1 < w_2 < \cdots < w_a$. Thus
\[
A_{\Pi(a,1^{n-a})}(t) = \sum_{w \in S_{n-a}}\binom{a+w_1-2}{a-1}
t^{\mathrm{des}(w)+ 1}.
\]

We give an explicit formula for $a = 2$.

\begin{prop}
\label{prop:multiEulerian} For $n \geq 3$, $A_{\Pi(2,1^{n-2})}(t)
= t\frac{d}{dt}A_{n-1}(t)$, where $A_{n-1}(t)$ is the usual
Eulerian polynomial. In other words,
\begin{equation}
\label{eq:des} \sum_{w \in S_{n-1}} w_1 t^{\mathrm{des}(w) + 1} =
\sum_{w \in S_{n-1}} (\mathrm{des}(w) + 1)t^{\mathrm{des}(w) + 1}.
\end{equation}
\end{prop}

We will prove this statement bijectively.  Let $w \in S_n$. A
\emph{circular descent} of $w$ is either a usual descent or the
index $n$ if $w_n
> w_1$. Define the \emph{circular descent number} $\cdes(w)$ as
the number of circular descents. Let $C_n$ denote the subgroup of
$S_n$ generated by the long cycle $c = (12\cdots n)$, written in
cycle notation. We have the following easy lemma.

\begin{lemma}
\label{lem:doublecosets} The statistic $\cdes$ is constant on
double cosets $C_n\backslash S_n/C_n$.
\end{lemma}
\begin{proof}
Left multiplication by the long cycle $c$ maps $w_1 w_2 \cdots
w_n$ to $(w_1+1)(w_2+1)\cdots(w_n+1)$ where ``$n+1$'' is
identified with ``1''.  Right multiplication by $c$ maps $w_1 w_2
\cdots w_n$ to $w_2w_3\cdots w_n w_1$.  The lemma is immediate
from the definition of circular descent number.
\end{proof}
Actually in the following we will use this Lemma for $S_{n+1}$.

\begin{proof}[Proof of Proposition~\ref{prop:multiEulerian}.]
If $w = w_1w_2\cdots w_{n} \in S_{n}$ and $w' = w_1 w_2 \cdots w_n
(n+1) \in S_{n+1}$ then $\cdes(w') = \des(w) + 1$.  By this and
our earlier comments, the left hand side of (\ref{eq:des}) is the
generating function for permutations in $S_{n+1}$ satisfying $w_1
< w_2$ and $w_{n+1} = n+1$, according to their circular descent
number.  Alternatively, we may view this as the $\cdes$-generating
function of right cosets $\tw \in S_{n+1}/C_{n+1}$ satisfying the
property that the two numbers $w_{i},w_{i+1}$ cyclically located
after $w_{i-1} = n+1$ satisfy $w_{i} < w_{i+1}$ for any
representative $w \in \tw$. Here the indices are taken modulo
$n+1$ and by Lemma~\ref{lem:doublecosets}, $\cdes(\tw) :=
\cdes(w)$ does not depend on the representative $w$ of $\tw$ and
so is well defined.

The right hand side of (\ref{eq:des}) is the generating function
for permutations $u = u_1u_2 \cdots u_n \in S_n$ satisfying $u_n =
n$ where one of the circular descents has been marked, again
according to circular descent number.  If $u_i
> u_{i+1}$ is the marked circular descent (where $i+1$ is to be taken modulo $n$),
then we can insert the number $n+1$ between $u_i$ and $u_{i+1}$ to
obtain a coset $\tv \in S_{n+1}/C_{n+1}$ with the same number of
circular descents.  If $v \in \tv$ and $v_i = n+1$ then we
automatically have $v_{i-1} > v_{i+1}$ (in fact this is a
$\cdes$-preserving bijection between cosets $\tv$ satisfying this
property and permutations of $u \in S_n$ satisfying $u_n = n$ with
a marked circular descent).  Let $c = (123\cdots(n+1))$ be the
generator of $C_{n+1}$ and consider $v' = c^{n + 1 -v_{i-1}} v$
for any $v \in \tv$ where $i$ is determined by $v_i = n+1$. Let
$\tv' = v'C$. By Lemma~\ref{lem:doublecosets}, $\cdes(\tv') =
\cdes(\tv)$. However, it is easy to see that $\tv'$ is exactly one
of the cosets which are enumerated by the left hand side of
(\ref{eq:des}).  Thus we obtain a $\cdes$-preserving bijection
$\tv \mapsto \tv'$ between two classes of cosets in
$S_{n+1}/C_{n+1}$, enumerated by the two sides of (\ref{eq:des}).

We illustrate the bijection with an example, where we will pick
representatives of appropriate cosets at our convenience.  Let $u
= 53162748 \in S_8$ satisfying $u_8 =8$ with marked circular
descent index $4$ corresponding to $u_4 = 6 > 2 = u_5$. This is an
object enumerated by the right hand side of (\ref{eq:des}), with
$n = 8$. Inserting ``$9$'' between ``$6$'' and ``$2$'' we obtain
$v = 531692748 \in S_9$. Multiplying on the left by $c^3$ (where
$c = (123456789)$) adds 3 to every value, changing ``6'' to ``9'',
 giving the permutation $v' = 864935172$.  Multiplying on the right
by $c^4$ we move the $9$ to the last position to get $w =
351728649 \in \tv' = v'C$, which satisfies $w_9 = 9$ and $w_1 <
w_2$.  This is exactly a permutation enumerated by the left hand
side of (\ref{eq:des}). Note that $\cdes(u) = 5 = \cdes(w)$ and
that all the steps can be reversed to give a bijection.
\end{proof}

It would be interesting if algebro-geometric proofs of some of our
results concerning (weighted) multi-Eulerian numbers could be
given; see Section~\ref{sec:grassmann}.

\section{Final Remarks} \label{sec:final_remarks}

\subsection{The $h$-vectors of alcoved polytopes}
%In Section~\ref{sec:f-vec}, we gave a formula for the $h$-vectors of
%the simplicial complex $\tilde{\G}_{2,n}$.
Let $\tilde{\G}_\P$ denote the simplicial complex associated to the
collection of simplices in the alcoved triangulation of $\P$.  The
triangulation $\G_\P$ of an alcoved polytope $\P$ is unimodular and
this implies that the $h$-polynomial $h(\tilde{\G}_\P,t)$ is a
non-negative polynomial.  We have given many interpretations for the
volume $h(\tilde{\G}_\P,1)$ of an alcoved polytope.  It would be
interesting to obtain statistics on these interpretations which give
$h(\tilde{\G}_\P,t)$.

In fact, the Erhart polynomial of $\P$ equals to the Hilbert
polynomial of the quotient ring $k[x_I]/\J_\P$ associated to the
polytope $\P$; see Section~\ref{sec:grobner_intro}. The
$h$-polynomial is the numerator of the associated generating
function.

\subsection{Relation with order polytopes}
Let $P$ be a poset naturally labelled with the numbers $[n]$.  The
order polytope $\oo_P \subset \R^n$ is defined by the inequalities
$0 \leq x_i \leq 1$ together with $x_i \leq x_j$ for every pair of
elements $i,j \in P$ satisfying $i > j$.  It is clear that $\oo_P$
is affinely equivalent to an alcoved polytope via the transformation
$y_i = x_i - x_{i-1}$.  The triangulation $\G_P$ of $\oo_P$ thus
obtained is known as the \emph{canonical triangulation} of $\oo_P$.

Denote by $\L_P$ the set of linear extensions of $P$;
see~\cite{EC1}.  The simplices of $\G_P$ can be labelled by $w \in
\L_P$ and this is compatible with the labelling used throughout
this paper.  It is known~\cite{Sta2} that we have
\[
h(\tilde{\G}_P,t) = \sum_{w \in \L_P} t^{\des(w)}.
\]
The descents counted by this $h$-vector, however, disagree with
the way we have been counting descents in this paper.  More
precisely, we have been concerned with the number of descents
$\des(w^{-1})$ of the inverse permutations labelling the
simplices.

\subsection{Weight polytopes for type $A_n$ and alcoved
polytopes} \label{sec:weight} A \emph{weight polytope} is the
convex hull of the weights which occur in some highest weight
representation of a Lie algebra.  For example,
the vertices of the hypersimplex $\H_{k,n}$ are exactly the weights
which occur in the $k^{th}$ fundamental representation of
$\mathrm{sl_n}$.

We will identify the integral weights $L$ of $\mathrm{sl_n}$ with
the integer vectors $(a_1,\ldots,a_n)$ satisfying $a_1 + \cdots +
a_n = l$, for some fixed $l$.  A weight polytope $\P_\lambda$ is
specified completely by giving the highest weight $\lambda =
(\lambda_1, \ldots, \lambda_n)$.  We will assume all the
coordinates $\lambda_i$ are non-negative.  A weight $\mu$ lies
inside $\P_\lambda$ if $\lambda$ dominates $\mu$ in the usual
sense: $\lambda_1 + \cdots + \lambda_i \geq \mu_1 + \cdots +\mu_i$
for all $i$.

Let $\mu, \nu \in L$ be two weights.  Define $U(\mu,\nu),
V(\mu,\nu) \in L$ by requiring that $U(I_\mu,I_\nu) =
I_{U(\mu,\nu)}$ and $V(I_\mu,I_\nu) = I_{V(\mu,\nu)}$, in the
notation of Section~\ref{sec:alcoved}.  Alternatively,
$U(\mu,\nu)$ and $V(\mu,\nu)$ are determined by requiring that
$\mu + \nu = U(\mu,\nu) + V(\mu,\nu)$ and $V(\mu,\nu) - U(\mu,\nu)
= \sum_i b_i \alpha_i$ for some $b_i \in \set{0,1}$, where
$\alpha_i = e_{i+1}-e_i$ are the simple roots.  We call
$\P_\lambda$ \emph{sort-closed} if the set of weights $\P_\lambda
\cap L$ is sort-closed, as in Theorem~\ref{thm:alcove}.

\begin{prop}
\label{prop:sorted_weight} A weight polytope $\P_\lambda$ is
sort-closed if and only if $\lambda$ is equal to $aw_i + bw_{i+1}$ for
non-negative integers $a,b$ and where $w_k =
(1,\ldots,1,0,\ldots,0)$ ($k$ 1's) are the fundamental weights.
\end{prop}

\begin{proof}
The `if' direction is easy as $\P_\lambda$ can be specified by $x_1
+ \cdots + x_n = i \cdot a +  (i+1) \cdot b$ and the inequalities $0
\leq x_j \leq a$ for each $j \in [1,n]$ and we can use
Theorem~\ref{thm:alcove}. For the other direction, we may assume the
highest weight $\lambda$ is a $n$-tuple with highest value $a$ and
lowest value $0$. Suppose $\lambda$ is not of the form of the
proposition, then there are two more values $b, c$ not equal to $a$
satisfying $a > b \geq c > 0$ so that $\lambda$ is of the form
$(a,\ldots,a,b,\ldots,b,c,\ldots,c,\ldots,0)$.

Explicitly construct a pair of weights $\delta =
(a,0,b,c,\ldots)$ and $\mu = (a-1,1,b+1,c-1,\ldots)$  where the
tails of the two $n$-tuples are identical, and $\delta$ is just a
permutation of the coordinates of $\lambda$. Both $\delta$ and
$\mu$ are dominated by $\lambda$ and hence lie in $\P_\lambda$.
However, $U(\delta,\mu) = (a,0,b+1,c-1,\ldots)$ does not lie in
$\P_\lambda$.
\end{proof}

Sturmfels~\cite[Chapter 14]{Stu} considered exactly this class of
sort-closed weight polytopes.
  The following
corollary follows immediately from Theorem~\ref{thm:alcove} and
Proposition~\ref{prop:sorted_weight}.

\begin{cor}
\label{cor:alcoved_weight} A weight polytope in the
$x$-coordinates with highest weight $\lambda$ is alcoved if and
only if $\lambda$ is of the form $a\omega_i + b\omega_{i+1}$ for
$a$ and $b$ non-negative integers. In particular, every weight
polytope for $A_2$ is alcoved.
\end{cor}

\subsection{Geometric motivation: degrees of torus
orbits} \label{sec:grassmann} Let $Gr_{k,n}$ denote the {\it
grassmannian manifold\/} of $k$-dimensional subspaces in the
complex linear space $\C^n$. Elements of $Gr_{k,n}$ can be
represented by $k\times n$-matrices of maximal rank $k$ modulo
left action of $GL_k$.  The $\binom{n}{k}$ maximal minors $p_I$ of
such a matrix, where $I$ runs over $k$-element subsets in
$\{1,\dots,n\}$, form projective coordinates on $Gr_{k,n}$, called
the {\it Pl\"ucker coordinates}. The map $Gr_{k,n}\to (p_I)$ gives
the {\it Pl\"ucker embedding\/} of the grassmannian into the
projective space $\mathbb{CP}^{\binom{n}{k} -1}$. Two points
$A,B\in Gr_{k,n}$ are in the same {\it matroid stratum\/} if
$p_I(A)=0$ is equivalent to $p_I(B)=0$, for all $I$.  The matroid
$\mm_A$ of $A$ has as set of bases $\set{I \in \binom{[n]}{k}
\mid p_I(A) \neq 0}$.

The complex torus $T=(\C\setminus\{0\})^n$ acts on $\C^n$ by
stretching the coordinates
$$
(t_1,\dots,t_n):(x_1,\dots,x_n)\mapsto(t_1 x_1,\dots,t_n x_n).
$$
This action lifts to an action of the torus $T$ on the grassmannian
$Gr_{k,n}$. This action was studied in~\cite{GGMS}. The authors
showed that the geometry of the closure of a torus orbit $X_A =
\overline{T\cdot A}$ depends (only) on the matroid stratum of $A$.
The variety $X_A$ is a toric variety and its associated polytope is
exactly the polytope $\P_{\mm_A}$ associated to the the matroid
$\mm_A$ from Section~\ref{sec:matroids}.  Our study of the volume of
the polytopes $\P_{\mm}$ was motivated by the well known fact
(see~\cite{Ful}) that
\[
\mathrm{deg}(X_A) = \mathrm{Vol}(\P_{\mm_A})
\]
where $\deg(X_A)$ denotes the degree of $X_A$ as a projective
subvariety of $\mathbb{CP}^{\binom{n}{k}-1}$ and ${\rm Vol}$
denotes the normalized volume with respect to the lattice
generated by the coordinate vectors $e_i$.  Note that by
definition only {\it representable} matroids $\mm$ arise as
$\mm_A$ in this manner.  Our Theorem~\ref{thm:matroid_tri} gives a
combinatorial description of the degree of a torus orbit closure
corresponding to a stratum of a sort-closed matroid.  In fact
Proposition~\ref{prop:grobner_multi} (and
Proposition~\ref{prop:sqfree}) shows that these torus orbit
closures are projectively normal, a fact known for all torus orbit
closures; see~\cite{Whi,Dab}.

It has been conjectured (see~\cite[Conjecture 13.19]{Stu}) that the
ideal of a smooth projectively normal toric variety is always
generated by quadratic binomials. The toric varieties associated to
simple alcoved polytopes give more examples of this.

\section{Appendix: Coherent triangulations and Gr\"{o}bner bases}
\label{sec:grobner_intro} We give a brief introduction to the
relationship between coherent triangulations of integer polytopes
and Gr\"{o}bner bases.  See Sturmfels~\cite{Stu} for further
details.

Let $k$ be a field and $k[x] = k[x_1,\ldots,x_n]$ be the
polynomial ring in $n$ variables. A total order $\prec$ on $\N^n$
is a \emph{term order} if it satisfies:
\begin{itemize}
\item
The zero vector is the unique minimal element.
\item
For any $\a,\,\b,\,\uc \in \N^n$, such that $\a \prec \b$ we have
$\a+\uc \prec \b+\uc$.
\end{itemize}
One way to create a term order is by giving a weight vector
$\w = (\w_1,\ldots,\w_n) \in \R^n$. Then for sufficiently generic
weight vectors a term order $\prec$ is given by $\b \prec \uc$ if
and only if $\w\cdot\b < \w\cdot\uc$\,.  In this situation we will say that
$\w$ represents $\prec$.

Given a polynomial $f \in k[x]$ one defines the initial monomial
${\rm in}_{\prec}(f)$ as the monomial $\underline{x}^\a$ with the
largest $\a$ under $\prec$. For an ideal $I$ of $k[x]$ one defines
the initial ideal ${\rm in}_\prec(I)$ as the ideal generated by
the initial monomials of elements of $I$.  The monomials which do
not lie in ${\rm in}_\prec(I)$ are called the standard monomials.
A finite subset $\g \subset I$ is a \emph{Gr\"{o}bner basis} for
$I$ with respect to $\prec$ if ${\rm in}_{\prec}(\g)$ generates
${\rm in}_\prec(I)$.  The Gr\"{o}bner basis is called reduced if
for two distinct elements $g, g' \in \g$, no term of $g'$ is
divisible by ${\rm in}_\prec(g)$.

Now let $\A = \set{\a_1,\a_2,\ldots,\a_n}$ be a finite subset of
$\Z^d$. We define a ring homomorphism $k[x] \to k[t_1^{\pm
1},\ldots,t_d^{\pm 1}]$ by
\[
x_k \longmapsto \underline{t}^{\a_k}.
\]
The kernel $\J_\A$ of this map is an ideal known as a \emph{toric
ideal}.

We now describe the relationship between Gr\"{o}bner bases of
$\J_\A$ and coherent triangulations of the convex hull of $\A$. For
any term order $\prec$, the initial complex $\D_\prec(\J_\A)$ of
$\J_\A$ is the simplicial complex defined as follows.  A subset $F
\subset \set{1,2,\ldots,n}$ is a face of $\D_\prec(\J_\A)$ if there
is no polynomial $f \in \J_\A$ such that the support of ${\rm
in}_\prec(f)$ is $F$. Thus the \emph{Stanley-Reisner ideal} of
$\D_\prec(\J_\A)$ is the radical of ${\rm in}_\prec(\J_\A)$.

A triangulation of a set $\A \in \Z^d$ (more specifically, its
convex hull) is \emph{coherent} if one can find a piecewise-linear
convex function $\nu$ on $\R^d$ such that the domains of linearity
are exactly the simplices of the triangulation.  Alternatively,
the triangulation is coherent if one can find a `height' vector
$\w$ such that the projection of the `lower' faces of the convex
hull of $\set{(\a_1,\w_1),(\a_2,\w_2),\ldots,(\a_n,\w_n)}$ is
exactly the triangulation.  We will denote such a triangulation of
$\A$ by $\D_\w(\A)$.  The function $\nu$ and the vector $\w$ can
be related by setting $\w_n = \nu(\a_n)$.

The main result we will need is the following~\cite[Chapter
8]{Stu}:

\begin{theorem}
\label{thm:regtri} The coherent triangulations of $\A$ are the
initial complexes of the toric ideal $\J_\A$. More precisely, if $\w
\in \R^n$ represents $\prec$ for $\J_\A$ then $\D_\prec(\J_\A) =
\D_\w(\A)$.
\end{theorem}

In the case when $\J_\A$ is a homogeneous toric ideal we can say
more.

\begin{prop}
\label{prop:sqfree} Let $\A$ be such that $\J_\A$ is a homogeneous
toric ideal.  Then the initial ideal ${\rm in}_\prec(\J_\A)$ is
square-free if and only if the corresponding regular triangulation
$\D_\prec$ of $\A$ is unimodular.  In that case, let $Y_{\A}$ be the
projective toric variety defined by the ideal $\J_\A$.  Then
$Y_{\A}$ is projectively normal and the Hilbert polynomial of
$Y_{\A}$ equals to the Erhart polynomial of the convex hull of $\A$.
\end{prop}
In this last case, the $\prec$-standard monomials correspond
exactly to the simplices of the triangulation.

\end{document}